\documentclass[11pt]{article}
\usepackage{amstext,amssymb,amsmath,amsbsy}
\usepackage{hyperref}
\usepackage{amscd} 
\usepackage{amsfonts,esint}
\usepackage{indentfirst}
\usepackage{verbatim}
\usepackage{amsmath}
\usepackage{amsthm}
\usepackage{enumerate}
\usepackage{graphicx}
\usepackage{subfig}
\usepackage{color}
\usepackage[OT1]{fontenc}
\usepackage[latin1]{inputenc}
\usepackage[english]{babel}
\usepackage{amssymb}
\makeatletter
\DeclareRobustCommand\widecheck[1]{
	{\mathpalette\@widecheck{#1}}}
\def\@widecheck#1#2{
		\setbox\z@\hbox{\m@th$#1#2$}
	\setbox\tw@\hbox{\m@th$#1%
		\widehat{%
			\vrule\@width\z@\@height\ht\z@
			\vrule\@height\z@\@width\wd\z@}$}%
	\dp\tw@-\ht\z@
	\@tempdima\ht\z@ \advance\@tempdima2\ht\tw@ \divide\@tempdima\thr@@
	\setbox\tw@\hbox{%
		\raise\@tempdima\hbox{\scalebox{1}[-1]{\lower\@tempdima\box
				\tw@}}}%
	{\ooalign{\box\tw@ \cr \box\z@}}}
\makeatother
\newtheorem{theorem}{Theorem}
\newtheorem{lemma}{Lemma}

\newtheorem{proposition}{Proposition}

\newtheorem{remark}{Remark}

\newtheorem{definition}{Definition}
\numberwithin{equation}{section}

\usepackage{dsfont}

\newcommand{\proofend}{\hfill $\Box$ }

\newcommand{\dsp}{\displaystyle}

\newcommand{\supp}{\operatorname{supp}}
\newcommand{\dist}{\operatorname{dist}}

\newcommand{\dive}{\operatorname{div}}

\newcommand{\eps}{\varepsilon}

\newcommand{\loc}{_{loc}}

\newcommand{\mC}{\mathbb{C}}

\newcommand{\mN}{\mathbb{N}}

\newcommand{\mR}{\mathbb{R}}

\newcommand{\cH}{{\cal H}}

\newcommand{\hv}{\hat v}
\newcommand{\tv}{\widetilde v}
\newcommand{\hg}{\hat g}
\newcommand{\hG}{\hat G}
\newcommand{\hw}{\hat w}

\title{Discreteness of interior transmission eigenvalues revisited}

\author{Hoai-Minh Nguyen\thanks{hoai-minh.nguyen@epfl.ch, \'Ecole Polytechnique F\'ed\'erale de Lausanne,  SB MATHAA CAMA, Station 8, CH-1015 Lausanne, Switzerland.}  \quad 
	Quoc-Hung Nguyen\thanks{quoc-hung.nguyen@epfl.ch, \'Ecole Polytechnique F\'ed\'erale de Lausanne,  SB MATHAA CAMA, Station 8, CH-1015 Lausanne, Switzerland.}}

\begin{document}

\maketitle 

\begin{abstract}
This paper is devoted to the discreteness of the transmission eigenvalue problems. 
It is known that this problem is not self-adjoint and a priori estimates are non-standard and do not hold in 
general.  Two approaches are used. The first one is based on the multiplier technique and the second one is based on the Fourier analysis.  The key point of the analysis is to establish the compactness and the uniqueness for Cauchy problems under various conditions. Using these approaches, we are able to rediscover quite a few known discreteness results in the literature and  obtain various new results for which  only the information near the boundary are required and there might be no contrast of the coefficients on the boundary.


\end{abstract}

\tableofcontents

\section{Introduction}	

Let $d \ge 2$ and $\Omega$ be a smooth bounded simply connected open subset of $\mR^d$ and  denote $\Gamma=\partial \Omega$.
 Let $A_1,A_2$ be two real symmetric matrix-valued functions and  $\Sigma_1,\Sigma_2$ be two bounded  positive functions all defined in $\Omega$ such that $A_1$ and $A_2$ are uniformly elliptic  and $\Sigma_1, \Sigma_2$ are bounded below by a positive constant in $\Omega$, i.e., for some constant $\Lambda \ge 1$,  and for $j=1, 2$, 
\begin{equation}\label{condi1}
 \Lambda^{-1} |\xi|^2\leq \langle  A_j(x) \xi, \xi \rangle \leq \Lambda |\xi|^2  \quad \mbox{ for all } \xi \in \mR^d, \mbox{ for a.e. }  x\in \Omega, 
\end{equation}
 and
 \begin{equation}\label{condi2}
 \Lambda^{-1}  \leq \Sigma_j(x)\leq \Lambda \mbox{ for a.e. } x \in \Omega.  
 \end{equation}
Here and in what follows  $\langle \cdot, \cdot \rangle$ denotes the Euclidean scalar product.  A complex number $\lambda$ is called an eigenvalue of the interior transmission eigenvalue (ITE) problem associated with the pairs $(A_1, \Sigma_1)$ and $(A_2, \Sigma_2)$ in $\Omega$ if there is a non-zero pair of functions $(u_1, u_2) \in [H^1(\Omega)]^2$ satisfying the system
 \begin{equation}\label{pro1a}  
 \left\{\begin{array}{lll}
 \dive(A_1 \nabla u_1) - \lambda\Sigma_1 u_1= 0 ~~&\text{ in}~\Omega, \\[6pt]
 \dive(A_2 \nabla u_2) - \lambda\Sigma_2 u_2= 0 ~~&\text{ in}~\Omega, 
 \end{array} \right. 
 \end{equation}
 and the boundary conditions
 \begin{equation}\label{pro1b}
  u_1 =u_2, \quad  A_1 \nabla u_1\cdot \nu = A_2 \nabla u_2\cdot \nu ~\text{ on }\Gamma, 
 \end{equation}
 where $\nu$ denotes the outward normal unit vector on $\Gamma$. Such a pair  $(u_1, u_2)$ is called an eigenfunction pair of \eqref{pro1a} and \eqref{pro1b}.

\medskip 
The ITE problem was  introduced in the middle of eighties by  Kirsch in \cite{Kir}  and Colton and Monk in \cite{ColtonMonk}. One of its interests comes from the connection between the density of the range of the far field operator with the injectivity of  ITE problem. The ITE problem is known to be not self-adjoint. More severely, a priori estimates are non-standard and do not hold in general. These create many difficulties for the investigation and a source of interesting problems. Some important directions in the study of the ITE problems are the  discreteness of the spectrum   \cite{Ca0, Coltonkirsch, ColPao, Faierman, Kir,LVainberg1,LVainberg2,Luc, RS, Sylvester}, the  completeness of the generalized  eigenfunctions \cite{Bla, Luc}, and the Weyl-laws for the spectrum  \cite{Faierman, Hi, LVainberg2, Luc}. The reader can find a  review on the ITE problem in \cite{Ca3}. 


\medskip 
This paper is devoted to the discreteness of the eigenvalues for the ITE problem.  The discreteness of the spectrum of  \eqref{pro1a}-\eqref{pro1b} was established by 
Rynne and Sleeman in \cite{RS} and 
Colton, Kirsch, Paivarinta in \cite{Coltonkirsch} in the case $A_1 = A_2 =I$, the identity matrix,  $\Sigma_1 - \Sigma_2 = \Sigma_1 -1 \ge c > 0$ in $\Omega$ for some constant $c$. Using the fact that $\Sigma_1 - \Sigma_2 > c > 0$, they can transform  \eqref{pro1a}-\eqref{pro1b}  into  an elliptic equation  of fourth order  and then derive the discreteness of the spectrum of the ITE problem from the new one. The discreteness of the spectrum of the ITE problem involving only the information near the boundary was  established quite recently. 
Using the T-coercivity method, which is related to the (Banach-Necas-Babuska) inf-sup condition, Bonnet-BenDhia, Chesnel, and  Haddar in \cite{Bo} obtained the discreteness of the spectrum of the ITE problem under the two assumptions $i)$ $(A_1, \Sigma_1) = (I, 1)$ or $(A_2, \Sigma_2) = (I, 1)$  and $ii)$ $A_1 - A_2  \ge c I$ and  $\Sigma_1 - \Sigma_2 \ge c$ in a neighborhood of $\Gamma$ for some positive constant $c$.  
In  \cite{Sylvester}, Sylvester showed that the discreteness takes place if $A_1 = A_2 = I$ in $\Omega$ and $\Sigma_1- \Sigma_2 = \Sigma_1 - 1 \ge c > 0$ in a neighbourhood of $\Gamma$ for some constant $c$ using the theory of  upper triangular compact operators and an apriori estimate for the ITE problem. His method also works in the case $A_1 = A_2$ smooth  in $\Omega$.  In \cite{LVainberg1}, 
Lakshtanov and Vainberg obtained the discreteness of ITE problem under the assumption that  $(A_1, \Sigma_1) = (I, 1)$, $(A_1, A_2, \Sigma_1, \Sigma_2)$ is smooth in $\Omega$ and  satisfies the  so-called parameter-elliptic conditions. In a related direction, Lakshtanov and Vainberg  in \cite{LVainberg2}  obtained the discreteness of the ITE problem for the case $(A_1, A_2, \Sigma_1, \Sigma_2) = (I, I, 1, \Sigma_2)$ in $\Omega$,  $\Sigma_2$ is smooth in $\bar \Omega$, and $\Sigma_2 = 1$ and  $\partial_\nu \Sigma_2 \neq 0$ on $\Gamma$.


\medskip 

In this paper, we follow the standard strategy used in the context of  Laplace operator  to investigate  the discreteness of the eigenvalues for the ITE problem.  To this end, we first establish the well-posedness of the following system, for some $\gamma_0 \in \mC$,  
 \begin{equation}\label{pro1a*} 
 \left\{ \begin{array}{cl}
 \dive(A_1 \nabla u_1) - \gamma_0\Sigma_1 u_1= g_1 &\text{ in}~\Omega, \\[6pt]
 \dive(A_2 \nabla u_2) - \gamma_0\Sigma_2 u_2= g_2 &\text{ in}~\Omega, \\[6pt]
  u_1 =u_2, \quad  A_1 \nabla u_1\cdot \nu = A_2 \nabla u_2\cdot \nu & \text{ on }\Gamma, 
 \end{array}
 \right.
 \end{equation}
where $(g_1,g_2)$ is a given pair of functions  in an  appropriate space. We then define the operator 
\begin{equation}\label{def-T-*}
T(f_1, f_2) = (u_1, u_2)  \mbox{ where $(u_1, u_2)$ is the unique solution of \eqref{pro1a*} with $(g_1, g_2) = (\Sigma_1 f_1, \Sigma_2 f_2)$}. 
\end{equation}
The discreteness  of the eigenvalues for the ITE problem can be now derived from the discreteness of the spectrum of $T$.  In this paper, the well-posedness and the compactness of the operator $T$, for some {\it appropriate} choice $\lambda_0 \in \mC$ are obtained via the multiplier technique  and  the Fourier analysis. The first approach has roots from the work of H-M. Nguyen in \cite{Ng-WP} where he studies the well-posedness of the Helmholtz equations with sign changing coefficients. Various ideas in this paper were introduced  there. 
The Fourier analysis is quite natural in this setting and the idea is to compute solutions of Cauchy problems in half space and derive the information from this. 

\medskip 
Throughout this paper
\begin{equation}\label{notatioin-delta}
d_\Gamma(x): = \mbox{dist} (x, \Gamma) := \inf_{y \in \Gamma} | x- y| \mbox{ for } x \in \mR^d,
\end{equation}
and  for two matrices $M_1,M_2$, we denote $M_1\geq M_2$ if  $M_1-M_2$ is a nonnegative matrix.

\medskip

Our first result of the discreteness of the eigenvalues for the ITE problem is 

\begin{theorem}\label{thm1}   Assume that for some  $0 \le \alpha < 2$ and for some positive constant $c$,  
\begin{equation}\label{cond-F-thm1}
			A_1  - A_2  \geq c d_\Gamma^\alpha I \quad \mbox{ and }  \quad  \Sigma_1 - \Sigma_2 \geq 0  \quad \mbox{ in a neighborhood of $\Gamma$}. 
\end{equation}
The spectrum of \eqref{pro1a}-\eqref{pro1b} is discrete. 
\end{theorem}

No regularity assumption on $(A_1,A_2,\Sigma_1,\Sigma_2)$ is imposed in Theorem \ref{thm1}.  The conditions in \eqref{cond-F-thm1} are only required near $\Gamma$ and neither upper bound of $A_1 - A_2$ nor the one of  $\Sigma_1 - \Sigma_2$ is required. The proof of Theorem~\ref{thm1} is based on the multiplier technique and 
given in Section~\ref{sect-goodsign}.  To our knowledge, Theorem~\ref{thm1} is new even in the case $\alpha = 0$. Applying Theorem~\ref{thm1} for $\alpha = 0$, one recovers and extends the discreteness results  obtained by Bonnet-BenDhia, Chesnel, and Haddar in \cite{Bo} mentioned.

\medskip 
When the second inequality  in \eqref{cond-F-thm1} is reversed and $\Sigma_2$ is large enough, the discreteness also holds. More precisely, we have

\begin{theorem}\label{thm3}
Assume that for some  $0 \le \alpha < 2$ and  for some positive constant $c$,  
\begin{align}\label{cond-F-thm3}
A_1 - A_2 \geq c  d_\Gamma^\alpha  I  \quad \mbox{ in a neighborhood of $\Gamma$}.  
\end{align}
Let $\Lambda_1 \ge 1$ be such that 
\begin{equation}\label{thm3-part1}
 \Lambda_1^{-1}I \leq A_j  \leq \Lambda_1 I \mbox{ for } j =1, 2, \quad \mbox{ and } \quad 
\Lambda_1^{-1} \le \Sigma_1 \le \Lambda_1 \quad  \mbox{ in }  \Omega. 
\end{equation}
For every $0< \Lambda_2 < 1$ there exists $K(\Lambda_2)  > 0$ depending only on $\Lambda_1$, $\Lambda_2$, $c$, $\alpha$, and the neighborhood of $\Gamma$ in \eqref{cond-F-thm3} such that if 
\begin{equation}\label{cond-K1}
\Lambda_2 K \le \Sigma_2  \leq K \mbox{ in }  \Omega,
\end{equation}
for some $K > K (\Lambda_2)$,  then  the spectrum of \eqref{pro1a}-\eqref{pro1b} is discrete. 
\end{theorem}

As far as we know, Theorem~\ref{thm3} is new even in the case $\alpha = 0$.  The proof of Theorem~\ref{thm3} is based on the multiplier technique and presented in Section~\ref{sect-badsign}. Another related result of Theorem~\ref{thm3}  assuming only a mild  condition on $(\Sigma_1, \Sigma_2)$ is given in Proposition~\ref{pro-A1A2}; however, one requires that $A_1 \ge A_2 \ge c d_{\Gamma}^\alpha I$ in the whole domain $\Omega$.  
 
\medskip 

We next deal with the case $A_1 = A_2$ in a neighborhood of $\Gamma$. 

\begin{theorem}\label{thm2}  Assume that for some $0 \le \beta < 2$ and   for some positive constants $c$,   
\begin{equation}\label{cond1-F-thm2}
A_1 = A_2 \quad \mbox{ and } \quad  \Sigma_1 - \Sigma_2 \ge c d_\Gamma^{\beta} \quad \mbox{in a neighborhood of $\Gamma$}.
\end{equation}
The spectrum of \eqref{pro1a}-\eqref{pro1b} is discrete. 
\end{theorem}

Applying Theorem~\ref{thm2} for the case $A_1 = A_2$ smooth in $\Omega$, and $\beta = 0$, one recovers  and extends the results mentioned of Sylvester's in \cite{Sylvester}  and Lakshtanov and Vainberg's in \cite{LVainberg2}. 
In comparison with their results, the novelty of Theorem~\ref{thm2} relies on the two facts: $i)$ only the information near the boundary on  $(A_1, A_2, \Sigma_1, \Sigma_2)$ is assumed and $ii)$ no regularity on the coefficients is required.

\medskip

In the case $(A_1, A_2, \Sigma_1, \Sigma_2)$ is continuous in a neighborhood of $\Gamma$, we prove the following result which  involves the complementing condition due to Agmon, Douglas, and Nirenberg in \cite{ADNII}:

\begin{theorem}\label{thm4} Assume that $A_1,\, A_2, \, \Sigma_1, \, \Sigma_2$ are continuous in a neighbourhood of $\Gamma$ and the following two conditions hold: 
\begin{enumerate}

\item[i)] For all $x \in \Gamma$, $A_1(x),  \, A_2(x)$ satisfy the complementing boundary condition with respect to  $\nu(x)$, 

\item[ii)]  For all $x \in \Gamma$, $\big\langle  A_1(x) \nu(x), \nu(x) \big\rangle \Sigma_1(x) \neq \big\langle  A_2(x) \nu(x), \nu(x) \big\rangle \Sigma_2(x)$. 
\end{enumerate}
The spectrum of \eqref{pro1a}-\eqref{pro1b} is discrete. 
\end{theorem}

The complementing condition corresponding to the ITE problem  has the following explicit characterization given in  \cite[Proposition 1]{Ng-WP} (see e.g. \cite[Definition 2]{Ng-WP} for the definition of the complementing condition for the ITE problem).  

\begin{proposition}\label{pro-complementing}  Let $e$ be a unit vector in $\mR^d$ and let $A_1$ and $A_2$ be two constant positive symmetric matrices. Then  $A_1$ and $A_2$ satisfy the complementing boundary condition with respect to $e$ if and only if 
\begin{equation}\label{cond-complementary}
\langle A_2 e, e \rangle \langle A_2 \xi, \xi \rangle  - \langle A_2 e,  \xi \rangle^2 \neq  \langle A_1 e, e \rangle \langle A_1 \xi, \xi \rangle  - \langle A_1 e, \xi \rangle^2 \quad \forall \, \xi \in P \setminus \{0 \}, 
\end{equation}
where
\begin{equation*}
{\cal P} := \big\{\xi \in \mR^d; \langle \xi, e \rangle = 0  \big\}. 
\end{equation*}
\end{proposition}

We next give some comments on the conditions $i)$ and $ii)$ of Theorem~\ref{thm4}:  

\begin{itemize}
\item It is shown in \cite[Proposition 1]{Ng-WP} that if $A_1 > A_2$ then $A_1$ and $A_2$ satisfy the complementing boundary condition for all unit vector $e$. 

\item It is clear from \eqref{cond-complementary} that if $d=2$ then \eqref{cond-complementary} is equivalent to the condition  $\det A_1 \neq \det A_2$. 

\item Assume that $A_1$ and $A_2$ are isotropic on $\Gamma$, i.e., $A_1 = a_1 I$ and $A_2 = a_2 I$ for some (positive) functions $a_1$, $a_2$ on $\Gamma$. Then $i)$ and $ii)$ are equivalent to the conditions $a_1 \neq a_2$ and $a_1 \Sigma_1 \neq a_2 \Sigma_2$ on $\Gamma$. 
\end{itemize}

The proof of Theorem~\ref{thm4} is via Fourier analysis and given in Section~\ref{sect-complementing}.  Applying Theorem~\ref{thm4}  and Proposition~\ref{pro-complementing}, one recovers and extends the result of   Lakshtanov and Vainberg's in \cite{LVainberg1}. 

\medskip 
The discreteness of the eigenvalues of a general elliptic system equipped complementing boundary condition was investigated by Faierman \cite{Faierman} and references therein. The conditions given in  \cite{Faierman} are not explicit and require that the coefficients are of class $C^1(\bar \Omega)$.  The proof in \cite{Faierman} involved the work of Agmon, Douglis, and Nirenberg   in \cite{ADNII}. 

\medskip 
It is worth noting that  ITE problem is  invariant after  a change of variables. In fact, let $F$ be  a diffeomorphism from $\Omega$ into itself \footnote{In this paper, this means that $F$ is bijective and $F, F^{-1} \in C^1 (\bar \Omega)$. } such that $F(x)=x$ on $\Gamma$.   
Set  $U_1 = u_1 \circ F^{-1}$ in $\Omega$. By the rule of a change of variables  (see e.g. \cite[Lemma 2]{Ng-Complementary}), 
$u_1 \in H^1(\Omega)$ is a solution to the equation 
$$
\dive(A_1 \nabla u_1) - \lambda \Sigma_1 u_1= 0 \mbox{ in } \Omega
$$
if and only if 
\begin{equation*} 
\dive(F_{*}A_1 \nabla U_1) - \lambda F_{*}\Sigma_1 U_1= 0 \mbox{ in } \Omega, 
\end{equation*}
where 
\begin{equation}\label{TO*}
F_{*}A_1(y) = \frac{D F  (x)  A_1(x) D F ^{T}(x)}{\det D F(x)} \quad  \mbox{  and  }  \quad  F_{*}\Sigma_1(y) = \frac{\Sigma_1(x)}{\det D F(x)} \mbox{ with } x =F^{-1}(y). 
\end{equation}
Moreover, 
\begin{equation}\label{pro1*'}  
U_1 = u_1 \quad \mbox{ and } \quad  F_*A_1 \nabla U_1 \cdot \nu = A_1 \nabla u_1 \cdot \nu \mbox{ on } \Gamma. 
\end{equation}
It follows that  $(u_1,u_2) \in [H^1(\Omega)]^2$ is a solution of \eqref{pro1a}-\eqref{pro1b}  if and only if $(U_1, U_2)=(u_1\circ F^{-1}, u_2) \in [H^1(\Omega)]^2$ is a solution of the system
\begin{equation}\label{pro1*'}  
\left\{\begin{array}{cll}
\dive(F_{*}A_1 \nabla U_1) - \lambda F_{*}\Sigma_1 U_1= 0 ~~&\text{ in}~\Omega \\[6pt]
\dive(A_2 \nabla U_2) - \lambda\Sigma_2 U_2= 0 ~~&\text{ in}~\Omega \\[6pt]
U_1 =U_2,~~~ F_{*}A_1 \nabla U_1\cdot \nu = A_2 \nabla U_2\cdot \nu &\text{ on }\Gamma.
\end{array}\right.
\end{equation}
Moreover, $(u_1,u_2) \not=0$ if and only if  $(U_1,U_2)\not= 0$. Using this observation, we can extend the previous results in which the conditions on $(A_1, \Sigma_1)$ and $(A_2, \Sigma_2)$ are involved to the case where the same conditions hold for $(F_*A_1, F_*\Sigma_1)$ and $(A_2, \Sigma_2)$ for some diffeomorphism $F$ verifying the condition $F(x) = x$ on $\Gamma$ 
\footnote{One can use two diffeomorphisms $F_1, \, F_2$ and require the corresponding conditions on $({F_1}_*A_1, {F_1}_*\Sigma_1)$ and $({F_2}_*A_2, {F_2}_*\Sigma_2)$ to obtain the discreteness of the ITE problem. However, the same conditions hold by using the diffeomorphisms $F_1\circ F_2^{-1}, \, I$.}.  We state here  two results following from Theorems~\ref{thm1} and \ref{thm2} in this direction as  an illustration.  A variant of Theorem~\ref{thm3} is left to the reader. Theorem~\ref{thm4} is already invariant under a change of variables. 
 
\begin{proposition}\label{pro-thm1-1}   Assume that for some diffeomorphisms $F$ from $\Omega$ into itself such that $F(x) = x$ on $\Gamma$, and for  $0 \le \alpha < 2$ and for some positive constant $c$,  
\begin{equation}\label{cond-F-pro-thm1}
			F_*A_1  - A_2  \geq c d_\Gamma^\alpha I \quad \mbox{ and }  \quad  F_*\Sigma_1 - \Sigma_2 \geq 0  \quad \mbox{ in a neighborhood of $\Gamma$}. 
\end{equation}
The spectrum of \eqref{pro1a}-\eqref{pro1b} is discrete. 
\end{proposition}

\begin{proposition}\label{pro-thm2-1}  Assume that for some diffeomorphism $F$ from $\Omega$ into itself such that $F(x) = x$ on $\Gamma$, 
\begin{equation}\label{cond1-F-pro-thm2}
F_*A_1 = A_2 \quad \mbox{ and } \quad  F_*\Sigma_1 - \Sigma_2 \ge c d_\Gamma^{\beta} \quad \mbox{in a neighborhood of $\Gamma$},
\end{equation}
for some $0 \le \beta < 2$ and   for some positive constants $c$.  The spectrum of \eqref{pro1a}-\eqref{pro1b} is discrete. 
\end{proposition}


Applying Propositions~\ref{pro-thm1-1} and \ref{pro-thm2-1}, one is able to obtain the discreteness of the spectrum of  the ITE problem even in the case $A_1 - A_2$ changes the sign in a neighborhood of $\Gamma$. 



\medskip 
We now describe briefly the ideas of the proof of Theorems~\ref{thm1}, \ref{thm3}, \ref{thm2}, and \ref{thm4}. The proof of Theorems~\ref{thm1}, \ref{thm3}, and \ref{thm2}  are based on the multiplier technique while the proof of Theorem~\ref{thm4} is based on the Fourier analysis. The key point is to establish the well-posedness and the compactness of $T$  given by \eqref{def-T-*} for some $\gamma_0 \in \mC$. Concerning Theorem~\ref{thm1}, the well-posedness of and the compactness of $T$ take place for  $\gamma_0 = \lambda_0$ for some large positive $\lambda_0$.   The existence of a solution of the system \eqref{pro1a*} is based on a priori estimates (Lemma~\ref{lem1}) for a Cauchy problem with its roots from \cite{Ng-WP} and via the limiting absorption process (Lemma~\ref{lem-elementary}). To require  the information only near the boundary $\Gamma$, $\lambda_0$ is chosen to be {\bf large} and a standard exponential decay estimate for elliptic equations is involved (Lemma~\ref{lem-decay}). Since the exponent $\alpha$  might be positive, some weighted spaces are involved and the solutions are not in $[H^1(\Omega)]^2$ as usual.  The proof of the uniqueness of \eqref{pro1a*} faces some issues due to the lack of the regularity of the solutions. To overcome this, we introduce the concept  of viscosity solution (Definition~\ref{def-WeakSolution1}), the terminology is inspired by from the one of systems of conservation laws. The compactness of $T$ follows from the condition $0 \le \alpha < 2$ in \eqref{cond-F-thm1}. Concerning Theorem~\ref{thm3}, the well-posedness of and the compactness of $T$ take place for  $\gamma_0 = i \lambda_0$ for some {\bf small} positive $\lambda_0$.  
This modification is necessary for the proof of the uniqueness. As in the proof of Theorem~\ref{thm1}, a priori estimates hold for \eqref{Sys2} (Lemma~\ref{lem1-v2}).
Nevertheless, the uniqueness of \eqref{Sys2} even for smooth solutions does not follow directly from the a priori estimates as in the proof of Theorem~\ref{thm1}. Additional arguments are required in this case (Lemma~\ref{lem-thm3}).  Beside these points,  the proof of Theorem~\ref{thm3} follows similarly as  the one of Theorem~\ref{thm1}. The proof of Theorem~\ref{thm2} is somehow in the spirit of the one of Theorem~\ref{thm1} but with the following  key difference. While the same arguments as used in the proof of Theorem~\ref{thm1} imply the uniqueness of $T$ with $\gamma_0 = \lambda_0$ for large $\lambda_0$, they do not imply the compactness of $T$ since  the first  condition in \eqref{cond-F-thm1} does not hold. To be able to deal with the situation in which $A_1 = A_2$ in a neighborhood of $\Gamma$, we  make some modifications on $T$. The idea is to take into account  the fact from the a priori estimates $u_1 - u_2 $  is more regular than  $u_1$ and $u_2$. After the modification of $T$, the proof of Theorem~\ref{thm2} is in the spirit of the one of Theorem~\ref{thm1} but the functional spaces used in this case is different and somehow more involved and the theory of compact analytic operator is used.  The proof of Theorem~\ref{thm4} is via the Fourier analysis. 
The approach is based on the computation of solutions of the Cauchy problems in half space via Fourier analysis. We then use  local  charts and the exponential decays of the solutions to deal with the general case.   The analysis is also in the spirit of the work Agmon, Douglis, and Nirenberg \cite{ADNII}.


\medskip 
The analysis in this paper is devoted to the study of the Cauchy problem \eqref{pro1a} and \eqref{pro1b}. There is a connection between the study of the Cauchy problem and the Helmholtz equations with sign changing coefficients modelling negative index materials. This connection was pointed out by H-M. Nguyen in \cite{Ng-WP}. Various ideas in this paper were introduced  there.
It is worth noting that resonance might appear in various interesting potential applications of negative index materials such as superlensing using complementary media \cite{Ng-Superlensing}, cloaking using complementary media \cite{Ng-Negative-cloaking, MinhLoc2}, cloaking a source via anomalous localized resonance, see e.g. \cite{Ng-CALR-CRAS, Ng-CALR, Ng-CALR-finite}, and cloaking an object via anomalous localized resonance in \cite{Ng-CALR-object}. 

\medskip 
The paper is organized as follows.  The proof of Theorems~\ref{thm1}, \ref{thm3}, \ref{thm2}, and \ref{thm4} are given in Sections~\ref{sect-goodsign}, \ref{sect-badsign}, \ref{sect-A1A2}, and  \ref{sect-complementing} respectively.

\section{On the case $A_1 \ge A_2$ and $\Sigma_1 \ge \Sigma_2$ in a neighborhood of $\Gamma$ - Proof of Theorem~\ref{thm1}} \label{sect-goodsign}

The section is devoted to the proof of Theorem~\ref{thm1}. To this end, as mentioned in the introduction, we first establish the well-posedness of the following system, for some $\lambda_0 > 0$,  
 \begin{equation}\label{Sys1} 
 \left\{ \begin{array}{cl}
 \dive(A_1 \nabla u_1) - \lambda_0\Sigma_1 u_1= g_1 &\text{ in}~\Omega, \\[6pt]
 \dive(A_2 \nabla u_2) - \lambda_0\Sigma_2 u_2= g_2 &\text{ in}~\Omega, \\[6pt]
  u_1 =u_2, \quad  A_1 \nabla u_1\cdot \nu = A_2 \nabla u_2\cdot \nu & \text{ on }\Gamma,
 \end{array}
 \right.
 \end{equation}
for a given pair $(g_1,g_2)$  in an  appropriate space. We then define the operator 
\begin{equation}\label{def-T}
T_1(f_1, f_2) = (u_1, u_2)  \mbox{ where $(u_1, u_2)$ is the unique solution of \eqref{Sys1} with $(g_1, g_2) = (\Sigma_1 f_1, \Sigma_2 f_2)$}
\end{equation}
and prove the compactness for it. 

The proof of the well-posedness of \eqref{Sys1} is as follows. 
The existence of a solution of \eqref{Sys1} is based on a priori estimates (Lemma~\ref{lem1}) for a Cauchy problem with roots from \cite{Ng-WP} and via a limiting absorption process (Lemma~\ref{lem-elementary}). To require only the information near the boundary, $\lambda_0$ is chosen to be large and a standard exponential decay estimate for elliptic equations is involved (Lemma~\ref{lem-decay}). The solutions are not in $[H^1(\Omega)]^2$ as usual and hence it is not clear whether the a priori estimates hold for the difference of  two solutions. 
To overcome this issue,  we  introduce  the concept of viscosity solutions (Definition~\ref{def-WeakSolution1}).  We now present the details of the proof. 
We start with  the following result (Lemma~\ref{lem1}) which  plays an important role in the proof of Theorem~\ref{thm1}. Lemma~\ref{lem1} is a more detailed version  of \cite[Lemma 9]{Ng-WP} and  its  proof follows closely  from there.

\begin{lemma} \label{lem1} Let $g = (g_1, g_2)  \in [L^2(\Omega)]^2$ and $h \in H^{-1/2}(\Gamma)$,  and $\lambda_0 \ge  0$,  and let   $v = (v_1, v_2)  \in [H^1(\Omega)]^2$  satisfy 
\begin{equation}\label{sys1-lem1} 
\left\{ \begin{array}{cl}
\dive(A_1 \nabla v_1) - \lambda_0\Sigma_1 v_1=  g_1 & \mbox{ in } \Omega, \\[6pt]
\dive(A_2 \nabla v_2) - \lambda_0\Sigma_2 v_2= g_2 & \mbox{ in } \Omega, \\[6pt]
v_1 = v_2, \quad  A_1 \nabla v_1\cdot \nu =  A_2 \nabla v_2\cdot \nu + h &  \text{ on }\Gamma.
\end{array}
\right.
\end{equation}
We have, with $w = v_1 - v_2$, 
\begin{equation}\label{part1-uniqueness}
\int_{\Omega} \langle A_1 \nabla w, \nabla w \rangle + \lambda_0 \Sigma_1 |w|^2 = -  \int_{\Omega} (g_1 - g_2) \bar w  +   \int_{\Omega} \lambda_0 (\Sigma_2 - \Sigma_1) v_2 \bar w +  \int_{\Omega} \langle [A_2 - A_1] \nabla v_2, \nabla w \rangle
\end{equation}
and 
\begin{multline}\label{lem1-part2-0}
\int_{\Omega}  \langle [A_1 - A_2] \nabla v_2, \nabla w \rangle  + \langle [A_1 - A_2] \nabla v_2,  \nabla v_2 \rangle   +  \lambda_0 (\Sigma_1 - \Sigma_2) |v_2|^2  \\[6pt]
 = \int_{\Omega}  g_2 \bar w  - (\bar g_1 - \bar g_2)  v_2  +  \lambda_0 (\Sigma_2 - \Sigma_1) v_2 \bar  w + \int_{\Gamma} \bar h v_2.  
\end{multline}
As a consequence of \eqref{part1-uniqueness} and  \eqref{lem1-part2-0}, we obtain 
\begin{equation}\label{lem1-part1}
\int_{\Omega} \langle A_1 \nabla w, \nabla w \rangle +  \lambda_0 \Sigma_1 |w|^2  \le 4 {\cal M}(v, g, h)  +   \int_{\Omega}  \langle [A_2 - A_1] \nabla v_2, \nabla w \rangle + \lambda_0 (\Sigma_2 - \Sigma_1) v_2 \bar w, 
\end{equation}
and 
\begin{multline}\label{lem1-part2}
\int_{\Omega} \langle [A_1 - A_2] \nabla v_2,  \nabla v_2  \rangle  +  \lambda_0 (\Sigma_1 - \Sigma_2) |v_2|^2  \\[6pt]
 \le 4 {\cal M}(v, g, h)  +  \int_{\Omega} \big| \langle [A_2 - A_1] \nabla v_2 \nabla w \rangle \big| +  \lambda_0 \big|(\Sigma_2 - \Sigma_1)   v_2 \bar  w \big|, 
\end{multline}
where
\begin{equation}\label{def-M}
{\cal M}(v, g, h) : = \|v \|_{L^2(\Omega)}\|g \|_{L^2(\Omega)} +  \|h\|_{H^{-1/2}(\Gamma)} \|v \|_{H^{1/2}(\Gamma)}. 
\end{equation}

\end{lemma}

Here and in what follows in this paper, for a complex number $z$,  $\bar z$ denotes its conjugate.  

\medskip 

\noindent{\bf Proof.} We derive from the definition of $w$ that $w = 0 $ on $\Gamma$ and 
\begin{equation}\label{part0-uniqueness}
\dive(A_1 \nabla w) - \lambda_0  \Sigma_1 w = g_1 - g_2 -\lambda_0 (\Sigma_2 - \Sigma_1) v_2  +  \dive ([A_2 - A_1] \nabla v_2) \mbox{ in } \Omega.
\end{equation}
Multiplying this equation by $\bar w$  and integrating on $\Omega$, we have
\begin{equation*}
\int_{\Omega} \langle A_1 \nabla w, \nabla w \rangle + \lambda_0 \Sigma_1 |w|^2 = -  \int_{\Omega} (g_1 - g_2) \bar w  +   \int_{\Omega} \lambda_0 (\Sigma_2 - \Sigma_1) v_2 \bar w +  \int_{\Omega} \langle [A_2 - A_1] \nabla v_2, \nabla w \rangle; 
\end{equation*}
which is \eqref{part1-uniqueness}. 
Multiplying the equation of $v_2$ by $\bar w$ and  integrating on $\Omega$,  we obtain
\begin{equation}\label{part0-uniqueness-1}
\int_{\Omega} - \langle A_2 \nabla v_2, \nabla w \rangle  -  \lambda_0 \Sigma_2 v_2 \bar w =  \int_{\Omega} g_2 \bar w. 
\end{equation}
It is clear that 
\begin{equation}\label{part1-uniqueness-1}
- A_2 \nabla w =  (A_1 - A_2) \nabla w  -  A_1 \nabla w + (A_2 - A_1) \nabla v_2 +  (A_1 - A_2) \nabla v_2 \mbox{ in } \Omega, 
\end{equation}
and,  by \eqref{part0-uniqueness}, 
\begin{equation}\label{part2-uniqueness-1}
\dive(A_1 \nabla w) - \dive ([A_2 - A_1] \nabla v_2)  =  \lambda_0  \Sigma_1 w +  g_1 - g_2 -\lambda_0 (\Sigma_2 - \Sigma_1) v_2 \mbox{ in } \Omega.
\end{equation}
Since 
$$
\big(A_1 \nabla w - (A_2 - A_1) \nabla v_2 \big) \cdot \nu = \big( A_1 \nabla v_1 - A_2 \nabla v_2 \big) \cdot \nu = h \mbox{ on } \Gamma,
$$
it follows from \eqref{part0-uniqueness-1}, \eqref{part1-uniqueness-1}, and \eqref{part2-uniqueness-1} that 
\begin{multline*}
\int_{\Omega} \langle [A_1 - A_2] \nabla v_2, \nabla w \rangle  + \langle [A_1 - A_2] \nabla v_2,  \nabla v_2 \rangle   +  \lambda_0 (\Sigma_1 - \Sigma_2) |v_2|^2  \\[6pt]
 = \int_{\Omega}  g_2 \bar w  - (\bar g_1 - \bar g_2)  v_2  +  \lambda_0 (\Sigma_2 - \Sigma_1) v_2 \bar  w + \int_{\Gamma} \bar h v_2; 
\end{multline*}
which is \eqref{lem1-part2-0}.  

\medskip 
Assertion~\eqref{lem1-part1} and  \eqref{lem1-part2} are  direct consequences of \eqref{part1-uniqueness} and \eqref{lem1-part2-0}. The proof is complete. \proofend

\medskip 
In what follows, we denote, for $s > 0$, 
$$
\Omega_s : = \big\{x \in \Omega;  \; d_\Gamma(x) < s \big\}. 
$$ 

\medskip 
The following exponential decay property for elliptic equations is useful for our analysis. 

\begin{lemma}\label{lem-decay} Let $\lambda> 1$,   $f \in L^2(\Omega)$, $A$ be a matrix-valued function,  and $\Sigma$ be real function defined in $\Omega$ such that, for some $\Lambda > 1$,  
\begin{equation}\label{cond1-lem2}
\Lambda^{-1} \leq A  \leq \Lambda \quad \mbox{ and } \quad \Lambda^{-1}  \leq \Sigma \leq \Lambda
\mbox{ in }   \Omega.
\end{equation}
Let  $u \in H^1_{\loc}(\Omega)$ be a solution to the equation $\dive (A \nabla u)  - \lambda \Sigma u = f$ in $\Omega$. 
 For all $s > 0$, there exist two positive constants $c_1$ and $c_2$, depending only on $\Lambda$, $s$, and $\Omega$, such that  
 \begin{equation}\label{state-lem2}
 \| u \|_{H^1(\Omega \setminus \Omega_s)} \le  c_1 \exp(- c_2 \sqrt{\lambda}) \| u \|_{L^2(\Omega_s)} + c_1 \| f\|_{L^2(\Omega)}. 
 \end{equation}
\end{lemma}

\noindent{\bf Proof.} Lemma~\ref{lem-decay} is quite standard and its proof presented here is in the spirit of the one of \cite[Theorem 2.2]{HJNg1}.  Let $U \in H^1_0(\Omega)$ be the unique solution of the equation $\dive (A \nabla U)  - \lambda \Sigma U = f$ in $\Omega$. Then $\| U\|_{H^1(\Omega)} \le  C \| f\|_{L^2(\Omega)}$. 
Here and in what follows in this proof $C$ denotes a positive constant depending only on $\Lambda$, $s$,  and $\Omega$.
By considering $u-U$, without loss of generality, one might assume that $f = 0$ in $\Omega$; this is assumed from now.  Fix $\varphi \in C^2(\Omega)$  such that $\varphi = c s$ in $\Omega \setminus \Omega_s$ and $\varphi = 0$ in $\Omega_{s/2}$, and $|\nabla \varphi| \le c$  in $\Omega$ where  $c$ is a small positive constant defined later (the smallness of $c$ depends only on $\Lambda$ and $\Omega$, it is independent of $s$). 
Set $\phi(x) =  e^{\sqrt{\lambda} \varphi(x)}$ and   $v (x)= u(x) \phi(x)$ for $x \in \Omega$. 
Since $\dive (A \nabla u)  - \lambda \Sigma u = 0$ in $\Omega$, it follows that 
\begin{equation}\label{eq-v11}
\dive(A \nabla v) - \lambda \Sigma v = \sqrt{\lambda} \dive ( v A  \nabla \varphi) + \sqrt{\lambda} A \nabla v \nabla \varphi - \lambda v A \nabla \varphi \nabla \varphi  \mbox{ in } \Omega. 
\end{equation}
Note that, by Cacciopoli's inequality, 
$$
\| u\|_{H^1(\Omega_{s/2} \setminus \Omega_{s/4})} \le C \sqrt{\lambda} \| u \|_{L^2(\Omega_s)}. 
$$ 
Hence there exists  $\tau \in (s/4, s/ 2)$ such that $ \| u \|_{H^1\big(\partial (\Omega \setminus \Omega_\tau) \big)} \le C \sqrt{\lambda} \| u\|_{L^2(\Omega_{s})}$.   Multiplying \eqref{eq-v11} by $\bar v$, integrating by parts in $\Omega \setminus \Omega_{\tau}$, and using \eqref{cond1-lem2} and the fact that $\varphi = 0$ in $\Omega_{s/2}$,  we have
\begin{equation}\label{est1-decay1}
\int_{\Omega \setminus \Omega_{\tau}} |\nabla v|^2 +  \lambda  |v|^2 \le C \Big( \int_{\Omega \setminus  \Omega_{\tau}}  \sqrt{\lambda} |v| |\nabla v| |\nabla \varphi| + \lambda |v|^2 |\nabla \varphi|^2 \Big) + C \lambda\| u\|_{L^2(\Omega_{s})}^2. 
\end{equation}
 By taking $c$ small enough, one can absorb the first term in the RHS of \eqref{est1-decay1} by the LHS, and the conclusion follows. \proofend

\medskip 
We next recall the following compactness result from  \cite[Lemma 7]{Ng-WP}  \footnote{In fact, \cite[Lemma 7]{Ng-WP} is stated for $(u_n) \subset H^1(\Omega)$, however the result also holds  for $(u_n) \subset H^1_{\loc}(\Omega)$ and the proof is almost unchanged.}. 
\begin{lemma} \label{lem-compact-1}  Let  $0\le  \alpha < 2$  and $(u_n) \subset H^1_{\loc}(\Omega)$. Assume that 
\begin{equation}\label{lem3-cond}
 \sup_{n} \int_{\Omega}  \Big( d_\Gamma^\alpha |\nabla u_n|^2  + |u_n|^2 \Big) \, dx < + \infty. 
\end{equation}
Then $(u_n)$ is relatively compact in $L^2(\Omega)$. 
\end{lemma}

Using Lemma~\ref{lem-compact-1}, we can prove 

\begin{lemma} \label{lem-decay3} Let $\sigma$ be a bounded real function defined in $\Omega$ such that $\int_{\Omega} \sigma \neq   0$. There exists a positive constant $C$ such that 
\begin{equation}\label{est-decay3}
\int_\Omega |u|^2\leq C\left( \left| \int_{\Omega} \sigma |u|^2 \right| +\int_\Omega d_{\Gamma}^\alpha |\nabla u|^2\right) \quad \forall \, u \in H^1_{\loc}(\Omega). 
\end{equation}
\end{lemma}

\noindent{\bf Proof.} Lemma~\ref{lem-decay3} can be derived from Lemma~\ref{lem-compact-1} by a contradiction argument as follows. Assume that \eqref{est-decay3} does not hold. There exists a sequence $(u_n) \subset H^1_{\loc}(\Omega)$ such that 
\begin{equation}\label{contradiction-decay3}
n \left( \left| \int_{\Omega} \sigma |u_n|^2 \right| +\int_\Omega d_{\Gamma}^\alpha |\nabla u_n|^2  \right) \le  \int_\Omega |u_n|^2 = 1. 
\end{equation}
By Lemma~\ref{lem-compact-1}, one may assume that $u_n \to u$ strongly in $L^2(\Omega)$ and almost everywhere in $\Omega$, and $u_n \rightharpoonup u$ weakly in $H^1_{\loc}(\Omega)$. We deduce from \eqref{contradiction-decay3} that 
\begin{equation*}
\left| \int_{\Omega} \sigma |u|^2\right| +\int_\Omega d_{\Gamma}^\alpha |\nabla u|^2  = 0 \quad \mbox{ and } \quad \int_\Omega |u|^2 = 1. 
\end{equation*}
The first identity implies that $u = 0$ in $\Omega$ since $\Omega$ is connected; this contradicts the second identity. Hence \eqref{est-decay3} holds. \proofend

\medskip 
Using Lemma~\ref{lem-decay}, we obtain 

\begin{lemma} \label{lem-stability}  Let $0 \le \alpha < 2$ and assume that, for some $\tau > 0$ and $c>0$, 
\begin{equation}\label{cond-F-thm1-*}
			A_1  - A_2  \geq c d_\Gamma^\alpha I \quad \mbox{ and }  \quad  \Sigma_1 - \Sigma_2 \geq 0 \mbox{ in } \Omega_\tau,  
\end{equation}
There exists $\Lambda_0 > 1$ such that if $\lambda_0 \ge  \Lambda_0$, and $v= (v_1, v_2)  \in [H^1_{\loc}(\Omega)]^2$  satisfy the system 
 \begin{equation}\label{sys1-cor1} 
\left\{ \begin{array}{cl}
\dive(A_1 \nabla v_1) - \lambda_0\Sigma_1 v_1=  g_1 & \mbox{ in } \Omega, \\[6pt]
\dive(A_2 \nabla v_2) - \lambda_0\Sigma_2 v_2= g_2 & \mbox{ in } \Omega,
\end{array}
\right.
\end{equation}
and  \eqref{lem1-part1} and \eqref{lem1-part2} hold for some  $g = (g_1, g_2) \in [L^2(\Omega)]^2$ and $h \in H^{-1/2}(\Gamma)$,   
then, with $w = v_1 - v_2$,  
\begin{equation}\label{est-lem-stability}
\int_{\Omega_\tau} \langle [A_1 - A_2] \nabla v_2, \nabla v_2 \rangle   + \int_{\Omega \setminus \Omega_{\tau}} |\nabla v_2|^2 
 + \int_{\Omega} |\nabla w|^2  + |v_2|^2  
\le C \Big( {\cal M}(v, g, h) + \| g\|_{L^2(\Omega)}^2 \Big), 
\end{equation}
for some positive constant  $C$ independent of $v$, $g$, and $h$ where ${\cal M}(v, g, h)$ is defined by \eqref{def-M}. In particular if $v \in [H^1(\Omega)]^2$ is a solution of \eqref{sys1-lem1} then \eqref{est-lem-stability} holds. 
\end{lemma}

\noindent{\bf Proof.} Since $A_1 - A_2 \ge 0$ and $\Sigma_1 - \Sigma_2 \ge 0$ in $\Omega_\tau$, we have,  for all $\gamma > 0$, 
\begin{equation}\label{Cauchy}
\gamma \int_{\Omega_\tau} \langle [A_1 - A_2] \nabla w, \nabla w \rangle + \frac{1}{\gamma} \int_{\Omega_\tau} \langle [A_1 - A_2] \nabla v_2,  \nabla v_2  \rangle  \ge 2 \int_{\Omega_\tau} \Big| \langle [A_1 - A_2] \nabla w, \nabla v_2 \rangle \Big|
\end{equation}
and 
\begin{equation}\label{Cauchy1}
\gamma \int_{\Omega_\tau}  (\Sigma_1 - \Sigma_2) |w|^2 + \frac{1}{\gamma} \int_{\Omega_\tau} (\Sigma_1 - \Sigma_2) |v_2|^2  \ge 2 \int_{\Omega_\tau}  (\Sigma_1 - \Sigma_2) |w v_2|. 
\end{equation}
Taking $\gamma > 1$ and close to 1, and adding \eqref{lem1-part1} and \eqref{lem1-part2}, we derive that 
\begin{multline}\label{lem1-part1-s}
\int_{\Omega_\tau} \langle \nabla w, \nabla w \rangle  +    \lambda_0   |w|^2 + \langle [A_1 - A_2] \nabla v_2,  \nabla v_2  \rangle   \\[6pt]
  \le C {\cal M}(v, g, h)  +  C \int_{\Omega \setminus \Omega_\tau}  |\nabla v_2| |\nabla w|  + \lambda_0 |v_2| |w| + |\nabla v_2|^2 + \lambda_0|v_2|^2. 
\end{multline}
Here and in what follows in this proof, $C$ denotes a positive constant independent of $v, \, g, \, h$ and $\lambda_0$. 

On the other hand, applying Lemma~\ref{lem-multiplier} below, we have
\begin{equation*}
\int_{\Omega} \lambda_0  d_\Gamma^{\alpha +2}  |v_2|^2 \le C \int_{\Omega} d_\Gamma^\alpha |\nabla v_2|^2 + C \| g_2\|_{L^2(\Omega)}^2. 
\end{equation*}
This implies, by applying Lemma~\ref{lem-decay3} with $\sigma = d_\Gamma^{\alpha+2}$, 
\begin{equation}\label{claim1-s}
\int_{\Omega} |v_2|^2 \le C \int_{\Omega} d_\Gamma^\alpha |\nabla v_2|^2 + C \| g_2\|_{L^2(\Omega)}^2. 
\end{equation}
The conclusion now follows from \eqref{lem1-part1-s} and \eqref{claim1-s} by applying Lemma~\ref{lem-decay} to  $v_2$ and taking $\Lambda_0$ large enough. \proofend
 
 \begin{remark} \label{rem-negative} \fontfamily{m} \selectfont Using  \eqref{claim1-s}, one can weaken the assumption $\Sigma_1 - \Sigma_2 \ge 0$ in \eqref{cond-F-thm1-*} by $\Sigma_1 - \Sigma_2 \ge - \hat c $ in a neighborhood of $\Gamma$ for some small positive constant $\hat c$. 
\end{remark}

In the proof of Lemma~\ref{lem-stability}, we used the following result. 
\begin{lemma}\label{lem-multiplier} 
 Let  $\lambda>1$, $\alpha \ge 0$, $f \in L^2(\Omega)$,  $A$ be a matrix-valued function,  and $\Sigma$ be real function defined in $\Omega$ such that 
\begin{equation*}
\Lambda^{-1} \leq A  \leq \Lambda \quad \mbox{ and } \quad \Lambda^{-1}  \leq \Sigma \leq \Lambda
\mbox{ in }   \Omega.
\end{equation*}
Let  $u \in H^1_{\loc}(\Omega)$ be a solution to the equation $\dive (A \nabla u)  - \lambda \Sigma u = f$ in $\Omega$. 
We have
\begin{equation}\label{multiplier1}
\int_{\Omega} \lambda  d_\Gamma^{\alpha +2}  |u|^2 \le C \int_{\Omega} d_\Gamma^\alpha |\nabla u|^2 + C \| f\|_{L^2(\Omega)}^2
\end{equation}
and
\begin{equation}\label{multiplier2}
\int_{\Omega} d_\Gamma^{\alpha + 2} |\nabla u|^2  \le C\int_{\Omega} \lambda  d_\Gamma^{\alpha}  |u|^2 +  C \| f\|_{L^2(\Omega)}^2,  
\end{equation}
for some positive constant $C$ independent of $\lambda$, $f$, and $u$.
\end{lemma}

\noindent{\bf Proof.} We only prove \eqref{multiplier1}; the proof of \eqref{multiplier2} follows similarly.  Set $\Gamma_s = \partial (\Omega \setminus \Omega_s)$ and 
$\varphi_s (x) = \dist(x, \Gamma_s)^{\alpha/ 2 + 1}$. 
Multiplying the equation of $u$ by $\bar u \varphi_s^2$ and integrating in $\Omega \setminus \Omega_s$, we have 
\begin{equation*}
\Big|\int_{\Omega \setminus \Omega_s} \lambda |u \varphi_s|^2 - \int_{\Omega \setminus \Omega_s} |\nabla u|^2  \varphi_s^2 \Big| \le C \int_{\Omega \setminus \Omega_s} |\nabla u| |\nabla \varphi_s| |u \varphi_s| + |f| |u \varphi_s^2|. 
\end{equation*}
Since $ |\nabla u| |\nabla \varphi_s| |u \varphi_s| \le  |\nabla \varphi_s|^2   |\nabla u|^2+ |\varphi_s|^2 |u|^2 \le C \dist^\alpha(x, \Gamma_s) |\nabla u|^2 +  |\varphi_s|^2 |u|^2 $, it follows that
\begin{equation*}
\int_{\Omega \setminus \Omega_s} \lambda |u \varphi_s|^2 \le C \int_{\Omega \setminus \Omega_s} \dist(x, \Gamma_s)^{\alpha}  |\nabla u|^2  + |f|^2. 
\end{equation*}
Letting $s \to 0$, applying the Fatou lemma to the LHS and the dominated convergence theorem to the RHS, one obtains \eqref{claim1-s}. \proofend

\begin{remark} \fontfamily{m} \selectfont The same conclusion  of Lemma~\ref{lem-multiplier}  holds if instead of imposing that $\dive (A \nabla u ) - \lambda \Sigma u = f$ in $\Omega$, one assumes $\dive (A \nabla u ) - i \lambda \Sigma u = f$ in $\Omega$. 

\end{remark}

\medskip 

The next basic result involving   the limiting absorption process is

\begin{lemma} \label{lem-elementary} Let $\lambda_0 \ge  0$, $0 < \delta < 1$,   $g = (g_1, g_2) \in [L^2(\Omega)]^2$. There exists a unique solution $v_\delta = (v_{1, \delta}, v_{2, \delta})  \in [H^1(\Omega)]^2$  of the system 
 \begin{equation}\label{sys1-lem-elementary} 
 \left\{ \begin{array}{cl}
 \dive \big((1 + i \delta ) A_1 \nabla v_{1, \delta} \big) - \lambda_0 \Sigma_1 v_{1, \delta} - i \delta v_{1, \delta}=  g_1 &\text{ in}~\Omega, \\[6pt]
  \dive\big((1 - i \delta )A_2 \nabla v_{2, \delta} \big) -\lambda_0 \Sigma_2 v_{2, \delta} + i \delta v_{2, \delta} = g_2 &\text{ in}~\Omega, \\[6pt]
 v_{1, \delta} = v_{2, \delta}, \quad (1 + i \delta) A_1 \nabla v_{1, \delta} \cdot \nu =  (1 - i \delta) A_2 \nabla v_{2, \delta}\cdot \nu&  \text{ on }\Gamma.
 \end{array}
 \right.
 \end{equation}
Moreover,  
\begin{equation}\label{lem6-1}
\|v_\delta \|_{H^1(\Omega)}^2 \le \frac{C}{\delta} \| g\|_{L^2(\Omega)}\| v_{\delta} \|_{L^2(\Omega)}, 
\end{equation}
for some positive constant $C$ independent of $g$ and $\delta$.  Consequently, 
\begin{equation}\label{lem6-2}
\|v_{\delta} \|_{H^1(\Omega)} \le \frac{C}{\delta} \| g \|_{L^2(\Omega)}. 
\end{equation}
\end{lemma}

\noindent{\bf Proof.} Set
\begin{equation}\label{def-X}
X:= \big\{\phi = (\phi_1,\phi_2)\in [H^1(\Omega)]^2: \; \phi_1-\phi_2  \in H^1_0(\Omega) \big\}.
\end{equation}
Then  $X$ is a Hilbert space equipped the scalar product: 
\begin{equation}\label{def-scalar-X}
\langle (\varphi_1,\varphi_2), (\phi_1,\phi_2)\rangle_{X}:= \langle \varphi_1, \phi_1\rangle_{H^1(\Omega)} + \langle \varphi_2, \phi_2\rangle_{H^1(\Omega)}.
\end{equation} 
Let $a$ be the bilinear form on $X \times X$ defined as follows, for $\varphi = (\varphi_1,\varphi_1), \phi = (\phi_1,\phi_2)\in X$, 
\begin{multline}\label{def-a}
a(\varphi, \phi)  = \int_\Omega  (1 + i \delta) \langle A_1 \nabla \varphi_1, \nabla \phi_1 \rangle + \int_\Omega  \lambda_0 \Sigma_1 \varphi_1 \overline{\phi}_1 +   \int_\Omega i \delta \varphi_1\overline{\phi}_1 \\[6pt]
- \left[\int_\Omega  (1 -i \delta )  \langle A_2 \nabla \varphi_2, \nabla \phi_2\rangle + \lambda_0 \int_\Omega \Sigma_2 \varphi_2 \overline{\phi}_2 -  \int_\Omega i \delta \varphi_2\overline{\phi}_2  \right],
\end{multline}
and let $b$ be the linear form on $X$ given by
\begin{equation}\label{def-b}
b(\phi) = \int_{\Omega} - g_1 \bar \phi_1 + g_2 \bar \phi_2.
\end{equation}
It is easy to check that $a$ is a continuous on $X \times X$ and $b$ is continuous on $X$; moreover, 
$$
|b(\phi)| \le \|(g_1, g_2) \|_{L^2(\Omega)} \| (\phi_1, \phi_2)\|_{L^2(\Omega)}. 
$$
On the other hand, by considering the imaginary part of $a$, we have
\begin{equation*}
\Im \big( a(\phi, \phi) \big) \geq C\delta \|\phi \|_{X}^2;  
\end{equation*}
which implies the coercivity of $a$. By Lax-Milgram's theorem, there exists a unique $\varphi \in X$ such that 
\begin{equation}\label{LM}
a(\varphi, \phi) = b(\phi) \mbox{ for all } \phi \in X. 
\end{equation}
One can check that $(v_{1, \delta}, v_{2, \delta}) = \varphi$ is a solution of \eqref{sys1-lem-elementary}; moreover,  
if $(v_{1, \delta}, v_{2, \delta}) \in X$ is a solution of \eqref{sys1-lem-elementary} then $\varphi = (v_{1, \delta}, v_{2, \delta})$ is a solution of \eqref{LM}. The proof is complete.  \proofend

\medskip 
Let $\tau > 0$ be such that \eqref{cond-F-thm1} holds in $\Omega_\tau$.  Define 
\begin{equation}\label{def-H}
\cH (\Omega): = \Big\{ (u_1, u_2) \in [H^1_{\loc} (\Omega) \cap L^2(\Omega)]^2;  \; u_1 - u_2 \in H^1_0(\Omega), \; \| (u_1, u_2)\|_{\cH(\Omega)} < + \infty \Big\}, 
\end{equation}
where, for $u = (u_1, u_2)$ and $v = (v_1, v_2)$ in $[H^1_{\loc}(\Omega)]^2$, 
\begin{equation}\label{scalar-H}
\langle u, v \rangle_{\cH(\Omega)} = \int_{\Omega} \nabla (u_1 - u_2) \nabla (\bar v_1 - \bar v_2) +  \int_{\Omega \setminus \Omega_\tau} \nabla u \nabla \bar v + 
\int_{\Omega_\tau}  \big\langle [A_1 - A_2] \nabla u, \nabla v \big\rangle  +  \int_{\Omega} u \bar v . 
\end{equation}
One can check that $\cH(\Omega)$ equipped the scalar product in \eqref{scalar-H} is a Hilbert space. 

\medskip 
We have 
\begin{definition}\label{def-weak-1} Let $(g_1, g_2) \in [L^2(\Omega)]^2$. A pair of functions $(v_1, v_2) \in \cH(\Omega)$ is called a weak solution of \eqref{Sys1} if 
\begin{equation}\label{sys-def-weak-1}
\left\{ \begin{array}{cl}
\dive(A_1 \nabla v_1) - \lambda_0\Sigma_1 v_1=  g_1 & \mbox{ in } \Omega, \\[6pt]
\dive(A_2 \nabla v_2) - \lambda_0\Sigma_2 v_2= g_2 & \mbox{ in } \Omega, \\[6pt]
(A_1 \nabla v_1-  A_2 \nabla v_2) \cdot \nu = 0  &  \text{ on }\Gamma
\end{array}
\right.
\end{equation}
\end{definition}

\begin{remark} \fontfamily{m} \selectfont
Since $\dive (A_1 \nabla v_1-  A_2 \nabla v_2) \in L^2 (\Omega)$ and $A_1 \nabla v_1-  A_2 \nabla v_2 = A_1 \nabla (v_1 - v_2) + (A_1 - A_2) \nabla v_2 \in L^2(\Omega)$ by \eqref{scalar-H}. The last identity on $\Gamma$ in \eqref{sys-def-weak-1} makes sense. 
\end{remark}

Using Lemmas~\ref{lem1},  \ref{lem-stability}, and  \ref{lem-elementary} one can construct a weak solution $(v_1, v_2) \in \cH(\Omega)$ of \eqref{Sys1}. More precisely, we have 

\begin{proposition}\label{pro-existence1} Assume \eqref{cond-F-thm1} and let $\delta \in (0, 1)$,  $g = (g_1, g_2) \in [L^2(\Omega)]^2$ and $v_\delta = (v_{1, \delta}, v_{2, \delta}) \in [H^1(\Omega)]^2$ be the unique solution of \eqref{sys1-lem-elementary}. There exists $\Lambda_0 > 1$ such that if $\lambda_0 > \Lambda_0$ then  
\begin{equation}\label{part1-pro-E1}
\|v_{\delta} \|_{\cH(\Omega)} \le C \|g\|_{L^2(\Omega)},
\end{equation}
for some positive constant  $C$ independent of $\delta$ and $g$. 
As a consequence, there exists a weak solution $v= (v_1, v_2) \in \cH(\Omega)$ of \eqref{Sys1} such that 
\begin{equation}\label{part2-pro-E1}
\|v\|_{\cH(\Omega)} \le C \| g \|_{L^2(\Omega)}. 
\end{equation}
\end{proposition}

\noindent{\bf Proof.} Applying Lemma~\ref{lem-elementary}, we have
\begin{equation} \label{est1-pro-E1}
\|v_{\delta} \|_{H^1(\Omega)} \le \frac{C}{\delta} \|g \|_{L^2(\Omega)}. 
\end{equation}
Here and in what follows in the proof, $C$ denotes a positive constant independent of $\delta$ and $g$. Rewriting the system of $v_\delta$, we get 
 \begin{equation}\label{sys1-pro-E1} 
 \left\{ \begin{array}{cl}
 \dive \big(A_1 \nabla v_{1, \delta} \big) - \lambda_0 \Sigma_1 v_{1, \delta} =  g_{1, \delta} &\text{ in}~\Omega, \\[6pt]
 \dive\big(A_2 \nabla v_{2, \delta} \big) -\lambda_0 \Sigma_2 v_{2, \delta}  =  g_{2, \delta} &\text{ in}~\Omega, \\[6pt]
A_1 \nabla v_{1, \delta} \cdot \nu -  A_2 \nabla v_{2, \delta}\cdot \nu = h_\delta&  \text{ on }\Gamma, 
\end{array} \right.
\end{equation}
where 
\begin{equation}\label{def-g1g2-1}
g_{1, \delta} =  \frac{1}{1 +  i \delta} \Big(g_1 + i \delta v_{1, \delta} - i \delta \lambda_0 \Sigma_1 v_{1, \delta}\Big),  \quad g_{2, \delta} = \frac{1}{1 -  i \delta} \Big(g_2 - i \delta v_{2, \delta} + i \delta \lambda_0 \Sigma_2 v_{2, \delta}\Big),  \\[6pt] 
\end{equation}
and 
\begin{equation*}
h_{\delta} =  - i \delta (A_1 \nabla v_{1, \delta} + A_2 \nabla v_{2, \delta} ) \cdot \nu. 
\end{equation*} 
Denote $g_\delta = (g_{1, \delta}, g_{2, \delta})$. 
For $\lambda_0 \ge \Lambda_0$ and for large $\Lambda_0$, we have, by Lemmas~\ref{lem1} and \ref{lem-stability} and \eqref{est1-pro-E1}, 
\begin{equation}\label{est2-pro-E1}
\| v_\delta \|_{\cH(\Omega)}^2  \le C \Big({\cal M} (v_\delta, g_\delta, h_\delta) + \| g_\delta\|_{L^2(\Omega)}^2 \Big) \le C \| g\|_{L^2(\Omega)}^2 + C  \| g \|_{L^2(\Omega)} \| v \|_{L^2(\Omega)}. 
\end{equation}
Here we used the fact that, by  \eqref{est1-pro-E1}, 
\begin{equation*}
\|g_\delta \|_{L^2(\Omega)}\le  C  \|  g \|_{L^2(\Omega)}
\end{equation*}
and in addition by the trace theory and \eqref{lem6-1}, 
\begin{align*}
\| h_\delta\|_{H^{-1/2}(\Gamma)} \|v_\delta \|_{H^{1/2}(\Gamma)} \le &  C \delta \Big(\| g_\delta\|_{L^2(\Omega)} + \| v_\delta\|_{H^1(\Omega)} \Big) \|v_\delta \|_{H^1(\Omega)} \\[6pt] 
\le &  C \Big( \|g \|_{L^2(\Omega)}^2 +  \|g\|_{L^2(\Omega)} \| v_\delta \|_{L^2(\Omega)} \Big).
\end{align*}
We derive from \eqref{est2-pro-E1} that 
\begin{equation}\label{est3-pro-E1}
\|v_{\delta} \|_{\cH(\Omega)}  \le C \| g \|_{L^2(\Omega)}; 
\end{equation}
which is \eqref{part1-pro-E1}. Assertion \eqref{part2-pro-E1} follows from \eqref{part1-pro-E1} by a standard compactness argument. The details are left to the reader. \proofend

\medskip 
We next discuss the uniqueness of the weak solutions. To obtain the uniqueness, one needs to show that $(v_1, v_2) = (0, 0)$ if $(v_1, v_2)$ 
is a weak solution of \eqref{sys-def-weak-1} with $g_1 = g_2 = 0$. A natural way for this is to obtain \eqref{lem1-part1} and \eqref{lem1-part2} with $g_1 = g_2 = 0$ and $h = 0$ for $(v_1, v_2)$ and then to apply Lemma~\ref{lem-stability}.  Note that estimate \eqref{lem1-part1} follows from \eqref{part1-uniqueness} and  estimate \eqref{lem1-part2} follows from \eqref{lem1-part2-0}. It is not clear that \eqref{part1-uniqueness} and \eqref{lem1-part2-0} hold for $(v_1, v_2) \in \cH(\Omega)$ which satisfies \eqref{sys-def-weak-1}. One way to remedy the situation is to consider only solutions of \eqref{sys-def-weak-1} which come from the limiting absorption principle process given in Lemma~\ref{lem-elementary}. Even making this restriction, the proof of the uniqueness is still not straightforward  since \eqref{part1-uniqueness} and \eqref{lem1-part2-0} are non-linear with respect to $(v_1, v_2)$.  To overcome this difficulty, we use an kind of relaxation argument which we now introduce. 
Consider $(\hv_1, \hv_2) \in \cH(\Omega)$ and $(\hg_1, \hg_2) \in [L^2(\Omega)]^2$ such that 
\begin{equation}\label{sys1-Uniqueness-0} 
\left\{ \begin{array}{cl}
\dive(A_1 \nabla \hv_1) - \lambda_0 \Sigma_1 \hv_1=  \hg_1 & \mbox{ in } \Omega, \\[6pt]
\dive(A_2 \nabla \hv_2) - \lambda_0\Sigma_2 \hv_2= \hg_2 & \mbox{ in } \Omega, \\[6pt]
(A_1 \nabla \hv_1  -  A_2 \nabla \hv_2) \cdot \nu = 0 &  \text{ on } \Gamma. 
\end{array}
\right.
\end{equation}
Let $(g_1, g_2) \in [L^2(\Omega)]^2$ and $(v_1, v_2) \in [H^1(\Omega)]^2$ be a solution of the system
\begin{equation*}
\left\{ \begin{array}{cl}
\dive(A_1 \nabla v_1) - \lambda_0\Sigma_1 v_1=  g_1 & \mbox{ in } \Omega, \\[6pt]
\dive(A_2 \nabla v_2) - \lambda_0\Sigma_2 v_2= g_2 & \mbox{ in } \Omega, \\[6pt]
v_1  = v_2  &  \text{ on }\Gamma. 
\end{array}
\right.
\end{equation*}
Set $w = v_1 - v_2$ and $\hat w = \hat v_1 - \hat v_2$. We claim that 
\begin{equation}\label{lem1-part2-0-0-1}
\int_{\Omega} \langle A_1 \nabla w, \nabla \hat w \rangle + \lambda_0 \Sigma_1 w \bar \hw = -  \int_{\Omega} (g_1 - g_2) \bar \hw  +   \int_{\Omega} \lambda_0 (\Sigma_2 - \Sigma_1) v_2 \bar \hw +  \int_{\Omega} \langle [A_2 - A_1] \nabla v_2, \nabla \hw \rangle
\end{equation}
and
\begin{multline}\label{lem1-part2-0-0}
\int_{\Omega}  \langle [A_1 - A_2] \nabla v_2, \nabla \hw \rangle  +  \langle [A_1 - A_2] \nabla v_2,  \nabla \hv_2 \rangle   +  \lambda_0 (\Sigma_1 - \Sigma_2) v_2 \bar \hv_2  \\[6pt]
 = \int_{\Omega}  g_2 \bar \hw  - (\bar \hg_1 - \bar \hg_2)  v_2  +  \lambda_0 (\Sigma_2 - \Sigma_1) v_2 \bar  \hw. 
\end{multline}
Indeed, we have 
\begin{equation*}
\dive(A_1 \nabla w) - \lambda_0  \Sigma_1 w = g_1 - g_2 -\lambda_0 (\Sigma_2 - \Sigma_1) v_2  +  \dive ([A_2 - A_1] \nabla v_2) \mbox{ in } \Omega.
\end{equation*}
Multiplying this equation by $\bar \hw$ and integrating in $\Omega$, we obtain \eqref{lem1-part2-0-0-1}. On the other hand,  multiplying the equation of $v_2$ by $\bar \hw$ and  integrating in $\Omega$,  we derive that 
\begin{equation}\label{part0-uniqueness-1-0}
\int_{\Omega} - \langle A_2 \nabla v_2, \nabla \hw \rangle  -  \lambda_0 \Sigma_2 v_2 \bar \hw =  \int_{\Omega} g_2 \bar \hw. 
\end{equation}
We have 
\begin{equation*}
- A_2 \nabla \hw =  (A_1 - A_2) \nabla \hw  -  A_1 \nabla \hw + (A_2 - A_1) \nabla \hv_2 +  (A_1 - A_2) \nabla \hv_2 \mbox{ in } \Omega
\end{equation*}
and, as in \eqref{part2-uniqueness-1},
\begin{equation*}
\dive(A_1 \nabla \hw) - \dive ([A_2 - A_1] \nabla \hv_2)  =  \lambda_0  \Sigma_1 \hw +  \hg_1 - \hg_2 -\lambda_0 (\Sigma_2 - \Sigma_1) \hv_2 \mbox{ in } \Omega.
\end{equation*}
Since 
$$
\big(A_1 \nabla \hw - (A_2 - A_1) \nabla \hv_2 \big) \cdot \nu = \big( A_1 \nabla \hv_1 - A_2 \nabla \hv_2 \big) \cdot \nu = 0 \mbox{ on } \Gamma, 
$$
it follows from \eqref{part0-uniqueness-1-0} that 
\begin{multline*}
\int_{\Omega}  \langle [A_1 - A_2] \nabla v_2, \nabla \hw \rangle  +  \langle [A_1 - A_2] \nabla v_2,  \nabla \hv_2 \rangle   +  \lambda_0 (\Sigma_1 - \Sigma_2) v_2 \bar \hv_2  \\[6pt]
 = \int_{\Omega}  g_2 \bar \hw  - (\bar \hg_1 - \bar \hg_2)  v_2  +  \lambda_0 (\Sigma_2 - \Sigma_1) v_2 \bar  \hw. 
\end{multline*}
which is \eqref{lem1-part2-0-0}.  


\medskip 
We are ready to introduce the notion of viscosity solutions of \eqref{Sys1}.

\begin{definition}\label{def-WeakSolution1} Let $(g_1, g_2) \in [L^2(\Omega)]^2$. A weak solution  $(v_1, v_2) \in \cH(\Omega)$ of \eqref{sys-def-weak-1}  
 is called a viscosity  solution if \eqref{lem1-part2-0-0-1} and  \eqref{lem1-part2-0-0} hold for any  $(\hv_1, \hv_2) \in \cH(\Omega)$ and $(\hg_1, \hg_2) \in [L^2(\Omega)]^2$ satisfying \eqref{sys1-Uniqueness-0}. 
\end{definition}


We are now in the position to state and prove the crucial result of this section.

\begin{proposition}\label{pro-thm1} Assume \eqref{cond-F-thm1}. There exists $\Lambda_0 > 1$ such that for $\lambda_0 \ge \Lambda_0$ and  for  $g = (g_1, g_2) \in [L^2(\Omega)]^2$, there exists a unique viscosity solution $v = (v_1, v_2) \in \cH(\Omega)$ of \eqref{Sys1}. 
Moreover, 
\begin{equation}\label{est1-p-thm1}
\| v\|_{\cH (\Omega)} \le C \| g \|_{L^2(\Omega)}, 
\end{equation}
for some positive constant $C$ independent of $g$. 
\end{proposition}

\noindent{\bf Proof.} Let $v_\delta = (v_{1, \delta}, v_{2, \delta}) \in [H^1(\Omega)]^2$ be the unique solution of \eqref{sys1-lem-elementary} and set $w_\delta = v_{1, \delta} - v_{2, \delta}$. By Proposition~\ref{pro-existence1}, we have, for large $\Lambda_0$, 
\begin{equation}\label{stability-pro-thm1}
\| v_\delta \|_{\cH(\Omega)} \le C \| g\|_{L^2(\Omega)}. 
\end{equation}
Let $v$ be the weak limit of $(v_{\delta_n})$ in $\cH(\Omega)$ for some sequence $(\delta_n) \to 0$ such that $w_{\delta_n} \rightharpoonup w = v_1 - v_2$ weakly in $H^1(\Omega)$. Then $v$ is a weak solution of \eqref{Sys1}. We prove that $v$ is a viscosity solution.  Define $g_{1, \delta}, \, g_{2, \delta}$ via \eqref{def-g1g2-1}. 
Then,  by  \eqref{lem1-part2-0-0-1} and \eqref{lem1-part2-0-0}, 
\begin{equation}\label{lem1-part2-0-0-1-pro}
\int_{\Omega} \langle A_1 \nabla w_\delta, \nabla \hat w \rangle + \lambda_0 \Sigma_1 w_\delta \bar \hw = -  \int_{\Omega} (g_{1, \delta} - g_{2, \delta}) \bar \hw  +   \int_{\Omega} \lambda_0 (\Sigma_2 - \Sigma_1) v_{2, \delta} \bar \hw +  \int_{\Omega} \langle [A_2 - A_1] \nabla v_{2, \delta}, \nabla \hw \rangle
\end{equation}
and 
\begin{multline}\label{lem1-part2-0-0-pro}
\int_{\Omega} \langle [A_1 - A_2] \nabla v_{2, \delta}, \nabla \hw \rangle  + \langle [A_1 - A_2] \nabla v_{2, \delta},  \nabla \hv_2 \rangle   +  \lambda_0 (\Sigma_1 - \Sigma_2) v_{2, \delta} \bar \hv_2  \\[6pt]
 = \int_{\Omega}  g_{2, \delta} \bar \hw  - (\bar \hg_1 - \bar \hg_2)  v_{2, \delta}  +  \lambda_0 (\Sigma_2 - \Sigma_1) v_{2, \delta} \bar  \hw. 
\end{multline}
By \eqref{stability-pro-thm1},  it follows from \eqref{lem6-1} of Lemma~\ref{lem-elementary} that 
\begin{equation}\label{pro-pro}
\| v_{\delta}\|_{H^1(\Omega)} \le \frac{C}{\sqrt{\delta}} \| g \|_{L^2(\Omega)}. 
\end{equation}
By letting $\delta \to 0$ in \eqref{lem1-part2-0-0-1-pro} and  \eqref{lem1-part2-0-0-pro}, and using \eqref{pro-pro},  one obtains \eqref{lem1-part2-0-0-1} and \eqref{lem1-part2-0-0}
for $v$. Hence $v$ is  a  viscosity solution. 
The uniqueness of $v$ follows from the definition of viscosity solutions as follows. Assume that $(v_1, v_2)$ and $(\tv_1, \tv_2)$ are two viscosity solutions. Set 
$$
(V_1, V_2) = (v_1 - \tv_1, v_2 - \tv_2) \mbox{ in } \Omega. 
$$
Then $ (V_1, V_2)$ is a viscosity solution corresponding to the pair of data $(0, 0)$. By choosing $(\hv_1, \hv_2) = (V_1, V_2)$, one obtains \eqref{part1-uniqueness} and \eqref{lem1-part2-0} where $(v_1, v_2) = (V_1, V_2)$ and $(g_1, g_2, h) = (0, 0, 0)$.  It follows from Lemma~\ref{lem1} that \eqref{lem1-part1} and \eqref{lem1-part2} holds with $(g_1, g_2, h) = (0, 0, 0)$. This implies $(V_1, V_2) = (0, 0)$ by Lemma~\ref{lem-stability} for large $\lambda_0$. The proof is complete. \proofend

\medskip

Fix $\lambda_0 >\Lambda_0$ where $\Lambda_0$ is the constant  in Proposition~\ref{pro-thm1}. Define 
\begin{eqnarray}\label{def-T1-11}
\begin{array}{rcccc}
T_1: & [L^2(\Omega)]^2 &\to&  [L^2(\Omega)]^2\\[6pt]
& (f_1, f_2)  &  \mapsto & (u_1, u_2), 
\end{array}
\end{eqnarray}
where $(u_1, u_2) \in \cH(\Omega) $ is the unique viscosity solution of \eqref{Sys1} with $(g_1, g_2) = (\Sigma_1 f_1, \Sigma_2 f_2)$. Using \eqref{est1-p-thm1} and applying  Lemma~\ref{lem-compact-1}, we derive that  $T_1$ is compact. 

\medskip 
We are ready to give the 

\medskip
\noindent{\bf Proof of Theorem~\ref{thm1}.}  By  the theory of compact operator see,  e.g., \cite{BrAnalyse1}, the spectrum of $T_1$ is discrete. It is clear that if $(u_1, u_2)$ is an eigenfunction pair of the ITE problem corresponding to the eigenvalue $\lambda$ then $(u_1, u_2)$ is eigenfunction pair of $T_1$ corresponding to the eigenvalue $\frac{1}{\lambda - \lambda_0}$. Hence, the spectrum of the ITE problem is discrete. \proofend

\begin{remark} \label{rem-com1}  \fontfamily{m} \selectfont The condition $0\leq \alpha <  2$ is necessary to guarantee  the compactness of $T_1$. 
\end{remark}

\begin{remark}  \fontfamily{m} \selectfont Theorem~\ref{thm1} also holds if the condition $\Sigma_1 - \Sigma_2 \ge 0$
in \eqref{cond-F-thm1} is replaced by $\Sigma_1 - \Sigma_2 \ge - \hat c $ for some small positive constant $\hat c$ since Lemma~\ref{lem-stability} holds in this case (see Remark~\ref{rem-negative}). 
\end{remark}

\section{On the case $A_1 \ge A_2$ and $\Sigma_1$ not greater than $\Sigma_2 $ in a neighborhood of $\Gamma$} \label{sect-badsign}

This section contains two sections. In the first one, we deal with the case $A_1 \ge A_2$ and $\Sigma_1 \le \Sigma_2$ in a neighborhood of $\Gamma$ and  give the proof of Theorem~\ref{thm3}. In the second one, we deal with the case $A_1 \ge A_2$ globally in $\Omega$. The main result in this section is Proposition~\ref{pro-A1A2}.

\subsection{On the case $A_1 \ge A_2$ and $\Sigma_1 \le \Sigma_2$ in a neighborhood of $\Gamma$-Proof of Theorem~\ref{thm3}}

The section is devoted to the proof of Theorem~\ref{thm3}. As mentioned in the introduction, we first establish the well-posedness of the following system, for some $\lambda_0 > 0$ (small),  
 \begin{equation}\label{Sys2} 
 \left\{ \begin{array}{cl}
 \dive(A_1 \nabla u_1) - i  \lambda_0\Sigma_1 u_1= g_1 &\text{ in}~\Omega, \\[6pt]
 \dive(A_2 \nabla u_2) - i \lambda_0\Sigma_2 u_2= g_2 &\text{ in}~\Omega, \\[6pt]
  u_1 =u_2, \quad  A_1 \nabla u_1\cdot \nu = A_2 \nabla u_2\cdot \nu & \text{ on }\Gamma,
 \end{array}
 \right.
 \end{equation}
for a given pair $(g_1,g_2)$  in an  appropriate space. We then define the operator 
\begin{equation}\label{def-T}
T_2(f_1, f_2) = (u_1, u_2)  \mbox{ where $(u_1, u_2)$ is the unique solution of \eqref{Sys2} with $(g_1, g_2) = (\Sigma_1 f_1, \Sigma_2 f_2)$}
\end{equation}
and prove the compactness for $T_2$.  System \eqref{Sys2} is slightly different from the one \eqref{Sys1} where the constant $i$ appears in front of $\lambda_0$. This modification is necessary for the proof of the uniqueness. As in the proof of Theorem~\ref{thm1}, similar a priori estimates hold for \eqref{Sys2} (Lemma~\ref{lem1-v2}).
Nevertheless, the uniqueness of \eqref{Sys2} even for smooth solutions does not follow directly from the a priori estimates as in the proof of Theorem~\ref{thm1}. Additional arguments are required and the condition on the largeness of  $\Sigma_2$ is involved (Lemma~\ref{lem-thm3}).  Beside this point, the proof of Theorem~\ref{thm3} is in the spirit of the one of Theorem~\ref{thm1}.

\medskip 
We now process the proof of Theorem~\ref{thm3}.  We first establish a variant of Lemma~\ref{lem1}.

\begin{lemma}\label{lem1-v2} Let $\lambda_0 >  0$, $g = (g_1, g_2) \in [L^2(\Omega)]^2$, and $h \in H^{-1/2}(\Gamma)$. Assume that   $v = (v_1, v_2)  \in [H^1(\Omega)]^2$  satisfy the system 
 \begin{equation}\label{sys1-lem1-v2} 
 \left\{ \begin{array}{cl}
 \dive(A_1 \nabla v_1) -  i \lambda_0\Sigma_1 v_1=  g_1 &\text{ in}~\Omega, \\[6pt]
 \dive(A_2 \nabla v_2) - i  \lambda_0\Sigma_2 v_2= g_2 &\text{ in}~\Omega, \\[6pt]
 v_1 = v_2, \quad  A_1 \nabla v_1\cdot \nu = A_2 \nabla v_2\cdot \nu + h &  \text{ on }\Gamma.
 \end{array}
 \right.
 \end{equation}
We have, with $w = v_2 - v_1$, 
\begin{equation}\label{part1-uniqueness-v2}
\int_{\Omega} \langle A_1 \nabla w, \nabla w \rangle + i  \lambda_0 \Sigma_1 |w|^2 = -  \int_{\Omega} (g_1 - g_2) \bar w  +   \int_{\Omega} i \lambda_0 (\Sigma_2 - \Sigma_1) v_2 \bar w +  \int_{\Omega} \langle [A_2 - A_1] \nabla v_2, \nabla w \rangle
\end{equation}
and
\begin{multline}\label{lem1-part2-0-v2}
\int_{\Omega} \langle [A_1 - A_2]  \nabla v_2, \nabla w \rangle  + \langle [A_1 - A_2] \nabla v_2,  \nabla v_2 \rangle   + i \lambda_0 (\Sigma_2 - \Sigma_1) |v_2|^2  \\[6pt]
 = \int_{\Omega}  g_2 \bar w  - (\bar g_1 - \bar g_2)  v_2  + i \lambda_0 (\Sigma_1 + \Sigma_2) v_2 \bar  w + \int_{\Gamma} \bar h v_2.
\end{multline}
As a consequence of \eqref{part1-uniqueness-v2} and \eqref{lem1-part2-0-v2}, we obtain 
\begin{equation}\label{lem1-v2-state1}
\lambda_0 \Big|  \int_{\Omega} (\Sigma_2 -  \Sigma_1) |v_2|^2 \Big| \le 4 {\cal M}(v, g, h) + \lambda_0 \int_{\Omega} \Sigma_1 |w|^2 + 2 \lambda_0 \int_{\Omega} \Sigma_1 |v_2| |w|
\end{equation}
and
\begin{multline} \label{lem1-v2-state2}
\int_\Omega \langle A_1 \nabla w,  \nabla w \rangle + \int_{\Omega} \langle [A_1 - A_2] \nabla v_2, \nabla v_2 \rangle \\ \le 
4 {\cal M} (v, g, h)  + 2 \int_{\Omega} \big|\langle [A_2 - A_1] \nabla v_2, \nabla w \rangle \big| + 2 \lambda_0 \int_{\Omega} \Sigma_2 |v_2| |w|,   
\end{multline} 
where ${\cal M}(v, g, h)$ is defined in \eqref{def-M}. 
\end{lemma}

\noindent{\bf Proof.} The proof is in the spirit of the one of Lemma~\ref{lem1}; nevertheless, different coefficients appear in the conclusion due to the effect of the constant $i$ in the front of $\lambda_0$. For the convenient of the reader, we present the proof.  From the definition of $w$, we derive that $w = 0 $ on $\Gamma$ and 
\begin{equation}\label{part0-uniqueness-v2}
\dive(A_1 \nabla w) - i \lambda_0  \Sigma_1 w = g_1 - g_2 -  i \lambda_0 (\Sigma_2 - \Sigma_1) v_2  +  \dive ([A_2 - A_1] \nabla v_2) \mbox{ in } \Omega.
\end{equation}
Multiplying this equation by $\bar w$  and integrating on $\Omega$, we have
\begin{equation}
\int_{\Omega} \langle A_1 \nabla w, \nabla w \rangle + i  \lambda_0 \Sigma_1 |w|^2 = -  \int_{\Omega} (g_1 - g_2) \bar w  +   \int_{\Omega} i \lambda_0 (\Sigma_2 - \Sigma_1) v_2 \bar w +  \int_{\Omega} \langle [A_2 - A_1] \nabla v_2, \nabla w \rangle;  
\end{equation}
which is \eqref{part1-uniqueness-v2}.  Multiplying the equation of $v_2$ by $\bar w$ and  integrating on $\Omega$,  we obtain
\begin{equation}\label{part0-uniqueness-1-v2}
\int_{\Omega} - \langle A_2 \nabla v_2, \nabla w \rangle  -  i  \lambda_0 \Sigma_2 v_2 \bar w =  \int_{\Omega} g_2 \bar w. 
\end{equation}
It is clear that 
\begin{equation}\label{part1-uniqueness-1-v2}
- A_2 \nabla w =  (A_1 - A_2) \nabla w  -  A_1 \nabla w + (A_2 - A_1) \nabla v_2 +  (A_1 - A_2) \nabla v_2 \mbox{ in } \Omega, 
\end{equation}
and,  by \eqref{part0-uniqueness-v2}, 
\begin{equation}\label{part2-uniqueness-1-v2}
\dive(A_1 \nabla w) - \dive ([A_2 - A_1] \nabla v_2)  =  i  \lambda_0  \Sigma_1 w +  g_1 - g_2 - i \lambda_0 (\Sigma_2 - \Sigma_1) v_2 \mbox{ in } \Omega. 
\end{equation}
Since 
$$
\big(A_1 \nabla w - (A_2 - A_1) \nabla v_2 \big) \cdot \nu = \big( A_1 \nabla v_1 - A_2 \nabla v_2 \big) \cdot \nu = h \mbox{ on } \Gamma, 
$$
it follows from \eqref{part0-uniqueness-1-v2}, \eqref{part1-uniqueness-1-v2}, and \eqref{part2-uniqueness-1-v2} that 
\begin{multline}\label{part3-uniqueness-v2}
\int_{\Omega} \langle [A_1 - A_2] \nabla v_2, \nabla w \rangle  + \langle [A_1 - A_2] \nabla v_2,  \nabla v_2 \rangle   + i \lambda_0 (\Sigma_2 - \Sigma_1) |v_2|^2  \\[6pt]
 = \int_{\Omega}  g_2 \bar w  - (\bar g_1 - \bar g_2)  v_2  + i \lambda_0 (\Sigma_1 + \Sigma_2) v_2 \bar  w + \int_{\Gamma} \bar h v_2; 
\end{multline}
which is \eqref{lem1-part2-0-v2}.  

\medskip 
Subtracting \eqref{part1-uniqueness-v2} from  \eqref{part3-uniqueness-v2} and considering the imaginary part  yields 
\begin{equation}\label{lem1-v2-part3}
\lambda_0 \Big| \int_{\Omega} (\Sigma_2 -  \Sigma_1) |v_2|^2 \Big| \le 4 {\cal M}(v, g, h) + \lambda_0 \int_{\Omega} \Sigma_1 |w|^2 + 2 \lambda_0 \int_{\Omega} \Sigma_1 |v_2| |w|; 
\end{equation}
which is \eqref{lem1-v2-state1}.  Adding \eqref{part1-uniqueness-v2} and  \eqref{part3-uniqueness-v2}  and considering the real part implies 
\begin{multline}\label{lem1-v2-part3}
\int_\Omega \langle A_1 \nabla w,  \nabla w \rangle + \int_{\Omega} \langle [A_1 - A_2] \nabla v_2, \nabla v_2 \rangle \\ \le 
4 {\cal M}(v, g, h) +  2 \lambda_0 \int_{\Omega}  \Sigma_2 |v_2| |w| + 2 \Big| \int_{\Omega}\langle [A_2 - A_1] \nabla v_2, \nabla w \rangle \Big|.
\end{multline} 
which is \eqref{lem1-v2-state2}. The proof is complete. \proofend



\medskip
Here is a variant of Lemma~\ref{lem-decay}.

\begin{lemma}\label{lem-decay-2} Let $\lambda> 1$,  $f \in L^2(\Omega)$, $A$ be a matrix-valued function,  and $\Sigma$ be real function defined in $\Omega$ such that
\begin{equation*}
\Lambda^{-1} \leq A   \leq \Lambda \quad \mbox{ and } \quad \Lambda^{-1}  \leq \Sigma \leq \Lambda
\mbox{ in }   \Omega.
\end{equation*}
Let $u \in H^1_{\loc}(\Omega)$ be a solution to the equation $\dive (A \nabla u)  - i \lambda \Sigma u = f$ in $\Omega$. 
For all $s > 0$, there exist two positive constants $c_1$ and $c_2$, depending only on $\Lambda$, $s$,  and $\Omega$, such that  
 \begin{equation*}
  \| u \|_{H^1(\Omega \setminus \Omega_s)} \le c_1 \exp(-c_2  \sqrt{\lambda}) \| u \|_{L^2(\Omega_s)} + c_1 \| f\|_{L^2(\Omega)}. 
 \end{equation*}
\end{lemma}

\noindent{\bf Proof.} The proof of Lemma~\ref{lem-decay-2} is similar to the one of Lemma~\ref{lem-decay}. The details are omitted. \proofend

\medskip 

The following result is a variant of Lemma~\ref{lem-stability} and plays  an important role in the proof of the uniqueness in Theorem~\ref{thm3}. 

\begin{lemma}\label{lem-thm3} Assume that for some $c, \,  \tau>0$ and  for some  $0 \le \alpha < 2$, 
\begin{align}\label{cond-F-thm3-lem}
A_1 - A_2 \geq c  d_\Gamma^\alpha  I  \mbox{ in } \Omega_\tau.  
\end{align}
Let $\Lambda_1 \ge 1$ be such that, 
\begin{equation}\label{thm3-part1-lem}
 \Lambda_1^{-1}I \leq A_j \leq \Lambda_1 I \mbox{ for } j =1, 2, \quad \mbox{ and } \quad 
\Lambda_1^{-1} \le \Sigma_1 \le \Lambda_1 \mbox{ in  }  \Omega. 
\end{equation} For every $0< \Lambda_2 < 1$ there exists $K_1  > 1$ depending only on $\Lambda_1$, $\Lambda_2$, $c$, $\alpha$, and $\tau$  such that if 
\begin{equation}\label{cond-K1-lem}
\Lambda_2 K \le \Sigma_2  \leq K \mbox{ in }  \Omega,
\end{equation}
for some $K > K_1$,  then there exists  $\lambda_0 > 0$, depending on $K$, such that  if  $g = (g_1, g_2) \in [L^2(\Omega)]^2$, $h \in H^{-1/2}(\Gamma)$,  and  $v= (v_1, v_2) \in \cH(\Omega)$ verify 
 \begin{equation}\label{sys1-thm3-lem} 
 \left\{ \begin{array}{cl}
 \dive(A_1 \nabla v_1) - i \lambda_0\Sigma_1 v_1=    g_1 &\text{ in}~\Omega, \\[6pt]
 \dive(A_2 \nabla v_2) - i \lambda_0\Sigma_2 v_2=  g_2 &\text{ in}~\Omega, 
  \end{array}
\right.
\end{equation}
and \eqref{lem1-v2-state1} and \eqref{lem1-v2-state2}, we have
\begin{equation}\label{est1-p-thm3}
\| v \|_{\cH (\Omega)}^2 \le C \Big( {\cal M}(v, g, h) +  \| g\|_{L^2(\Omega)}^2 \Big), 
\end{equation}
for some positive constant $C$ independent of $v$,  $g $, and $h$ where ${\cal M} = {\cal M}(v, g, h)$ is defined in \eqref{def-M}.  In particular if $v \in [H^1(\Omega)]^2$ is a solution of \eqref{sys1-lem1-v2} then \eqref{est1-p-thm3} holds.  
\end{lemma}

\noindent{\bf Proof.} Fix $0< \eps < 1$ a small constant which is defined later and set $\lambda_0 = \eps K^{-1/2}$. Note that 
$\lambda_0 K = \eps K^{1/2}$ is large if $K$ is large; this fact is assumed from now.  We derive from \eqref{lem1-v2-state1} and Lemma~\ref{lem-decay-2} that, for large $K$,   
\begin{equation*}
\int_{\Omega} \lambda_0 K |v_{2}|^2 \le C  \int_{\Omega} \lambda_0 |v_{2} | |w| + \lambda_0 |w|^2 +  C {\cal M}(v, g, h). 
\end{equation*}
Here and in what follows in the proof, $C$ denotes a positive constant depending only on $\Lambda_1$, $\Lambda_2$,  $\Omega$,  and $\tau$; it is independent of $K$ and $\eps$. 
This implies, for large $K$,  
\begin{equation}\label{pro-thm3-p2-1}
\int_{\Omega} \lambda_0 K |v_{2}|^2 \le C  \int_{\Omega}  \lambda_0 |w|^2 +  C {\cal M}(v, g, h). 
\end{equation}
Using the fact, for $\gamma > 0$,  
\begin{equation*}
\gamma \int_{\Omega_\tau} \langle [A_1 - A_2] \nabla w, \nabla w \rangle + \frac{1}{\gamma} \int_{\Omega_\tau} \langle [A_1 - A_2] \nabla v_2,  \nabla v_2  \rangle  \ge 2 \int_{\Omega_\tau} \Big| \langle [A_1 - A_2] \nabla w, \nabla v_2 \rangle \Big|
\end{equation*}
and taking $\gamma > 1$ and close to 1, we derive from \eqref{lem1-v2-state2}  that 
\begin{equation}\label{pro-thm3-p2-1-1}
\int_\Omega \langle A_1 \nabla w,  \nabla w \rangle + \int_{\Omega_\tau} \langle [A_{1} - A_{2}] \nabla v_{2} , \nabla v_{2} \rangle \le  C  \int_{\Omega} \lambda_0 K |v_{2}| |w| + C\int_{\Omega \setminus \Omega_\tau} |\nabla v_2|^2 + C {\cal M}(v, g, h).  
\end{equation}
Since 
$$
2\lambda_0 K |v_{2}| |w| \le \lambda_0 K^{3/2} |v_2|^2 + \lambda_0 K^{1/2} |w|^2, 
$$
it follows from  \eqref{pro-thm3-p2-1} and \eqref{pro-thm3-p2-1-1} that 
\begin{multline}\label{pro-thm3-p3}
\int_\Omega \langle A_1 \nabla w,  \nabla w \rangle + \int_{\Omega_\tau} \langle [A_{1} - A_{2}] \nabla v_{2} , \nabla v_{2} \rangle \\[6pt] 
\le  C  \int_{\Omega} \lambda_0 K^{1/2}  |w|^2 + C\int_{\Omega \setminus \Omega_\tau} |\nabla v_2|^2 +  C K^{1/2}{\cal M}(v, g, h).  
\end{multline}
Choosing $\eps$ small enough, one can absorb the first  term in the RHS of \eqref{pro-thm3-p3} by the first term of the LHS of \eqref{pro-thm3-p3} (recall that $\lambda_0 K^{1/2} = \eps$) and obtains 
\begin{equation}\label{pro-thm3-p4-111}
\int_\Omega \langle A_{1} \nabla w,  \nabla w \rangle + \int_{\Omega_\tau} \langle [A_{1} - A_{2}] \nabla v_{2} , \nabla v_{2} \rangle
\le C\int_{\Omega \setminus \Omega_\tau} |\nabla v_2|^2 + C K^{1/2}{\cal M}(v, g, h).  
\end{equation}
Similar to \eqref{claim1-s}, we have
\begin{equation}\label{claim1-s-*}
\int_{\Omega} |v_2|^2  \le C \int_{\Omega} d_\Gamma^\alpha |\nabla v_2|^2 + C \| g_2\|_{L^2(\Omega)}^2. 
\end{equation}
Estimate~\eqref{est1-p-thm3} now follows from \eqref{pro-thm3-p4-111} and \eqref{claim1-s-*} by applying Lemma~\ref{lem-decay-2} with $v = v_2$,  $A = A_2$, $\Sigma = \Sigma_2/ K$, and $\lambda = \lambda_0 K = \eps K^{1/2}$.  
\proofend

\medskip 
The following result  is a variant of Lemma~\ref{lem-elementary} and its proof uses essentially Lemma~\ref{lem-thm3}. 

\begin{lemma} \label{lem-elementary-v2} 
Assume that for some $c, \, \tau>0$ and for some $0 \le \alpha < 2$, 
\begin{align}\label{cond-F-thm3-lem}
A_1 - A_2 \geq c  d_\Gamma^\alpha  I  \mbox{ in } \Omega_\tau.  
\end{align}
Let $\Lambda_1 \ge 1$ be such that
\begin{equation}\label{thm3-part1-lem}
 \Lambda_1^{-1}I \leq A_j \leq \Lambda_1 I \mbox{ for } j =1, 2, \quad \mbox{ and } \quad 
\Lambda_1^{-1} \le \Sigma_1 \le \Lambda_1 \mbox{ in  }  \Omega. 
\end{equation} For every $0< \Lambda_2 < 1$ there exists $K_1  > 1$ depending only on $\Lambda_1$, $\Lambda_2$, $c$, $\alpha$, and $\tau$  such that if 
\begin{equation}\label{cond-K1-lem}
\Lambda_2 K \le \Sigma_2  \leq K \mbox{ in }  \Omega,
\end{equation}
for some $K > K_1$,  then there exists  $\lambda_0 > 0$, depending on $K$, 
 such that for  $ \delta \in (0, 1)$,   $g = (g_1, g_2) \in [L^2(\Omega)]^2$, there exists a unique solution $v_\delta = (v_{1, \delta}, v_{2, \delta})  \in [H^1(\Omega)]^2$  of the system 
 \begin{equation}\label{sys1-lem-elementary-v2} 
 \left\{ \begin{array}{cl}
 \dive \big((1 +  \delta ) A_1 \nabla v_{1, \delta} \big) - i \lambda_0 \Sigma_1 v_{1, \delta} =  g_1 &\text{ in}~\Omega, \\[6pt]
  \dive\big( A_2 \nabla v_{2, \delta} \big) - i \lambda_0 \Sigma_2 v_{2, \delta}= g_2 &\text{ in}~\Omega, \\[6pt]
 v_{1, \delta} = v_{2, \delta}, \quad (1 +  \delta) A_1 \nabla v_{1, \delta} \cdot \nu =  A_2 \nabla v_{2, \delta}\cdot \nu&  \text{ on }\Gamma.
 \end{array}
 \right.
 \end{equation}
Moreover,  
\begin{equation}\label{estimate-11}
\|v_{\delta} \|_{H^1(\Omega)}^2 \le \frac{C}{\delta} \Big(\| g \|_{L^2(\Omega)}\| v_{\delta} \|_{L^2(\Omega)} +  \| g \|_{L^2(\Omega)}^2 \Big), 
\end{equation}
for some positive constant $C$ independent of $g$ and $\delta$.  Consequently, 
\begin{equation*}
\|  v_{\delta} \|_{H^1(\Omega)} \le \frac{C}{\delta} \| g \|_{L^2(\Omega)}. 
\end{equation*}
\end{lemma}

\noindent{\bf Proof.} Applying Lemma~\ref{lem-thm3} with $A_1  = (1+\delta) A_1$ and $A_2 = A_2$, there exists $K_1  > 1$ depending only on $\Lambda_1$, $\Lambda_2$, $c$, $\alpha$, and $\tau$  such that if 
\begin{equation}\label{cond-K1-lem-lem}
\Lambda_2 K \le \Sigma_2  \leq K \mbox{ in }  \Omega,
\end{equation}
for some $K > K_1$,  there exists  $\lambda_0 > 0$, depending on $K$, such that if $v_\delta = (v_{1, \delta}, v_{2, \delta}) \in [H^1(\Omega)]^2$ ($0< \delta < 1$) is a solution of \eqref{sys1-lem-elementary-v2} then 
\begin{multline}
\|v_{\delta} \|_{H^1(\Omega \setminus \Omega_{\tau/2})}^2 +  \| v_{1, \delta} - v_{2, \delta}\|_{H^1 (\Omega)}^2 + \|v_{\delta} \|_{L^2(\Omega)}^2 + 
\int_{\Omega_\tau}  \big\langle \big[(1 + \delta)A_1 - A_2 \big] \nabla v_{ \delta}, \nabla v_{ \delta} \big\rangle  \\[6pt]
\le C {\cal M}(v_\delta, g, 0) + C \| g\|_{L^2(\Omega)}^2,  
\end{multline}
for some positive constant $C$ independent of $g, v_\delta$, and $\delta$. This implies the uniqueness of $v_\delta$ and estimate \eqref{estimate-11}. 
The existence of $v_\delta$ follows from Fredholm's theory and can be  proceeded as follows. Define
\begin{eqnarray}
\begin{array}{rcccc}
T_{1, \eps}: & [L^2(\Omega)]^2 &\to&  [L^2(\Omega)]^2\\[6pt]
& (g_1, g_2)  &  \mapsto & (v_1, v_2), 
\end{array}
\end{eqnarray}
where $(v_1, v_2) \in [H^1(\Omega)]^2 $ is the unique solution of \eqref{Sys1} where $(A_1, A_2)$ and $(\Sigma_1, \Sigma_2)$  are replaced by 
$\big((1 + \delta)A_1, A_2 \big)$ and $(2 \eps^{-1}, \eps^{-1})$ for small $\eps$.  $T_{1, \eps}$ is well-defined since $T_1$ given in \eqref{def-T1-11} is well-defined ($\alpha = 0$ in this case). 
The existence of $v_\delta$ follows from the uniqueness of $v_\delta$ by applying the Fredholm theory for the operator $I - T_{1, \eps} \circ B: [L^2(\Omega)]^2 \to [L^2(\Omega)]^2$ where $B: [L^2(\Omega)]^2 \to [L^2(\Omega)]^2$ is defined by 
$$
B (v_1, v_2) = (i \lambda_0 \Sigma_1 v_1 - 2 \eps^{-1}  v_1, i \lambda_0 \Sigma_2 v_2 -   \eps^{-1} v_2).
$$
Clearly $B$ is invertible if $\eps$ is small enough. Note that  $I - T_{1, \eps} \circ B$  is injective by the uniqueness of \eqref{sys1-lem-elementary-v2}. By the Fredholm theory  it is bijective. Hence for any $g \in [L^2(\Omega)]^2$ there exists $u \in [L^2(\Omega)]^2$ such that $u - T_{1, \eps} \circ B (u) = B^{-1} g$. Set $v = u - B^{-1} g$. Then $T_{1, \eps} \circ B (v + B^{-1} g) =  v$. In other words,  $v$ is a solution of \eqref{sys1-lem-elementary-v2}.  The proof is complete. \proofend

\medskip 
We are now in the position to give the definition of a weak solution of \eqref{Sys1}. 
Let $\tau > 0$ be such that \eqref{cond-F-thm3} holds in $\Omega_\tau$.  Define $\cH(\Omega)$ as in \eqref{def-H} with the scalar product given in \eqref{scalar-H}. The notion of weak solution  is similar to the one in Definition~\ref{def-weak-1} and  is given in 
\begin{definition}\label{def-weak-2} Let $(g_1, g_2) \in [L^2(\Omega)]^2$. A pair of functions $(v_1, v_2) \in \cH(\Omega)$ is called a weak solution of \eqref{Sys1} if 
\begin{equation}\label{sys-def-weak-2}
\left\{ \begin{array}{cl}
\dive(A_1 \nabla v_1) - i \lambda_0\Sigma_1 v_1=  g_1 & \mbox{ in } \Omega, \\[6pt]
\dive(A_2 \nabla v_2) - i \lambda_0\Sigma_2 v_2= g_2 & \mbox{ in } \Omega, \\[6pt]
(A_1 \nabla v_1-  A_2 \nabla v_2) \cdot \nu = 0  &  \text{ on }\Gamma.
\end{array}
\right.
\end{equation}
\end{definition}


Using  Lemmas~\ref{lem-thm3} and \ref{lem-elementary-v2}, one can construct a weak solution $(v_1, v_2) \in \cH(\Omega)$ of \eqref{Sys2}. More precisely, we have 

\begin{proposition}\label{pro-existence1-v2} Let $\delta \in (0, 1)$,  $g = (g_1, g_2) \in [L^2(\Omega)]^2$,  and $v_\delta = (v_{1, \delta}, v_{2, \delta}) \in [H^1(\Omega)]^2$ be the unique solution of \eqref{sys1-lem-elementary-v2}. Assume 
\eqref{cond-F-thm3-lem}, \eqref{thm3-part1-lem}, and \eqref{cond-K1-lem}. Let $K > K_1 > 1$ and  $\lambda_0 > 0$ be as in Lemma~\ref{lem-elementary-v2}.  We have
\begin{equation}\label{part1-pro-E1-v2}
\| v_{\delta} \|_{\cH(\Omega)} \le C \| g \|_{L^2(\Omega)}, 
\end{equation}
for some positive constant  $C$  independent of $\delta$ and $(g_1, g_2)$.  As a consequence, there exists a weak solution $v = (v_1, v_2) \in \cH(\Omega)$ of \eqref{Sys2} such that 
\begin{equation}\label{part2-pro-E1-v2}
\|v \|_{\cH(\Omega)} \le C \| g \|_{L^2(\Omega)}. 
\end{equation}
\end{proposition}

\noindent{\bf Proof.} The proof is similar to the one of Proposition~\ref{pro-existence1}. However, instead of using Lemmas~\ref{lem1}, \ref{lem-stability}, and \ref{lem-elementary}, one applies Lemmas~\ref{lem1-v2},  \ref{lem-thm3}, and \ref{lem-elementary-v2}. The details are left to the reader. \proofend

\medskip As in the proof of Theorem~\ref{thm1} in Section~\ref{sect-goodsign}, we  introduce the concept of viscosity solutions using a relaxation argument. Let $(\hv_1, \hv_2) \in \cH(\Omega)$ and $(\hg_1, \hg_2) \in [L^2(\Omega)]^2$ be such that 
\begin{equation}\label{sys1-Uniqueness-0-v2} 
\left\{ \begin{array}{cl}
\dive(A_1 \nabla \hv_1) - i \lambda_0 \Sigma_1 \hv_1=  \hg_1 & \mbox{ in } \Omega, \\[6pt]
\dive(A_2 \nabla \hv_2) - i \lambda_0\Sigma_2 \hv_2= \hg_2 & \mbox{ in } \Omega, \\[6pt]
(A_1 \nabla \hv_1 -   A_2 \nabla \hv_2) \cdot \nu &  \text{ on }\Gamma. 
\end{array}
\right.
\end{equation}
Let $(g_1, g_2) \in [L^2(\Omega)]^2$ and $(v_1, v_2) \in [H^1(\Omega)]^2$ be a weak solution of the system
\begin{equation*}
\left\{ \begin{array}{cl}
\dive(A_1 \nabla v_1) - \lambda_0\Sigma_1 v_1=  g_1 & \mbox{ in } \Omega, \\[6pt]
\dive(A_2 \nabla v_2) - \lambda_0\Sigma_2 v_2= g_2 & \mbox{ in } \Omega, \\[6pt]
v_1  = v_2  &  \text{ on }\Gamma. 
\end{array}
\right.
\end{equation*}
Set $w = v_1 - v_2$ and $\hat w = \hat v_1 - \hat v_2$.  Involving the same arguments  used to derive \eqref{lem1-part2-0-0-1} and \eqref{lem1-part2-0-0}, we have   
\begin{equation}\label{lem1-part2-0-0-v2-1}
\int_{\Omega} \langle A_1 \nabla w, \nabla \hw \rangle + i  \lambda_0 \Sigma_1 w \bar \hw = -  \int_{\Omega} (g_1 - g_2) \bar \hw  +   \int_{\Omega} i \lambda_0 (\Sigma_2 - \Sigma_1) v_2 \bar \hw +  \int_{\Omega} \langle [A_2 - A_1] \nabla v_2, \nabla \hw \rangle
\end{equation}
and
\begin{multline}\label{lem1-part2-0-0-v2}
\int_{\Omega} \langle (A_1 - A_2) \nabla v_2, \nabla \hw \rangle  + \langle (A_1 - A_2) \nabla v_2,  \nabla \hv_2 \rangle   + i   \lambda_0 (\Sigma_2 - \Sigma_1) v_2 \bar \hv_2  \\[6pt]
 = \int_{\Omega}  g_2 \bar \hw  - (\bar \hg_1 - \bar \hg_2)  v_2  + i \lambda_0 (\Sigma_1 + \Sigma_2) v_2 \bar  \hw. 
\end{multline}

\medskip 
We are ready now to define the notion of viscosity solution of \eqref{Sys2}.

\begin{definition}\label{def-WeakSolution1-v2} Let $(g_1, g_2) \in [L^2(\Omega)]^2$. A weak solution  $(v_1, v_2) \in \cH(\Omega)$ of \eqref{Sys2}  
 is called a viscosity  solution if  \eqref{lem1-part2-0-0-v2-1} and  \eqref{lem1-part2-0-0-v2} hold for any  $(\hv_1, \hv_2) \in \cH(\Omega)$ and $(\hg_1, \hg_2) \in [L^2(\Omega)]^2$ which satisfy  \eqref{sys1-Uniqueness-0-v2}. 
\end{definition}

Applying Lemma~\ref{lem1-v2} and using Lemma~\ref{lem-elementary-v2}, we can  prove

\begin{proposition}\label{pro-thm3} Assume 
\eqref{cond-F-thm3-lem}, \eqref{thm3-part1-lem}, and \eqref{cond-K1-lem}. Let $K_1 > 1$ and $\lambda_0 > 1$ be as in Lemma~\ref{lem-elementary-v2}. For  $g = (g_1, g_2) \in [L^2(\Omega)]^2$, there exists a unique viscosity solution $v = (v_1, v_2) \in \cH(\Omega)$ of \eqref{Sys2}. 
Moreover, 
\begin{equation}
\| v\|_{\cH (\Omega)} \le C \| g \|_{L^2(\Omega)}, 
\end{equation}
for some positive constant $C$ independent of $(g_1, g_2)$. 
\end{proposition}

\noindent{\bf Proof.} The proof of Proposition~\ref{pro-thm3} is similar to the one of Proposition~\ref{pro-thm1}. The details are left to the reader. 
\proofend

\medskip 

Fix $\lambda_0$ in Proposition~\ref{pro-thm3}. Define 
\begin{eqnarray}
\begin{array}{rcccc}
T_2: & [L^2(\Omega)]^2 &\to&  [L^2(\Omega)]^2 \\[6pt]
& (f_1, f_2)  &  \mapsto & (u_1, u_2), 
\end{array}
\end{eqnarray}
where $(u_1, u_2) \in \cH(\Omega)$ is the unique viscosity solution of \eqref{Sys2}. The compactness of $T_2$ is  a consequence of Lemma~\ref{lem-compact-1}. 

\medskip 
We are ready to give the 

\medskip
\noindent{\bf Proof of Theorem~\ref{thm3}.}   By  the theory of compact operator, the spectrum of $T_2$ is discrete. It is clear that if $(u_1, u_2)$ is an eigenfunction pair of the ITE problem corresponding to $\lambda$ then $(u_1, u_2)$ is eigenfunction pair of $T$ corresponding to the eigenvalue $\frac{1}{\lambda - i\lambda_0}$. Hence, the spectrum of the ITE problem is discrete. \proofend

\begin{remark} \label{rem-com2}  \fontfamily{m} \selectfont As in the proof of Theorem~\ref{thm1},  the condition $0\leq \alpha <  2$ in Theorem~\ref{thm3} is necessary to guarantee  the compactness of $T_2$. 
\end{remark}

\subsection{On the case $A_1 \ge A_2$ globally in $\Omega$}

The main result of this section is 
\begin{proposition}\label{pro-A1A2} 	Assume that 
	\begin{equation}\label{L2.1-cond1*}
A_1 - A_2 \geq c d_{\Gamma}^\alpha I ~~\text{in }~\Omega  \quad \text{ and } \quad \int_{\Omega}(\Sigma_1-\Sigma_2)\not =0,
	\end{equation}
 for some $0\leq\alpha< 2$ and for some $c>0$. The spectrum of \eqref{pro1a}-\eqref{pro1b} is discrete. 
\end{proposition}

In Proposition~\ref{pro-A1A2} the condition   $A_1 - A_2 \geq c d_{\Gamma}^\alpha I$ is required  in the whole  domain $\Omega$ not in a neighborhood of $\Gamma$; this is different from the context of  Theorem~\ref{thm3}. Applying Proposition \ref{pro-A1A2} with $\alpha = 0$ and $(A_1, \Sigma_1) = (I, 1)$ in $\Omega$, one rediscovers the discreteness result given in \cite[Theorem 4.4]{Bo}. 

\medskip 
The proof of Proposition~\ref{pro-A1A2} is based on the following lemma

\begin{lemma} \label{lem-stability-*}   Assume that for some $c>0$ and for some $0 \le \alpha < 2$, 
\begin{equation}\label{cond-F-thm1-*-pro}
			A_1  - A_2 \geq c d_\Gamma^\alpha \mbox{ in } \Omega  \quad \mbox{ and }  \quad  \int_{\Omega} (\Sigma_1 - \Sigma_2) \neq 0. 
\end{equation}
There exists $\Lambda_0 >  0$ such that if $0 < \lambda_0 <   \Lambda_0$, and $v = (v_1, v_2)  \in [H^1_{\loc}(\Omega)]^2$  verifies 
 \begin{equation}\label{sys1-cor1} 
\left\{ \begin{array}{cl}
\dive(A_1 \nabla v_1) - i \lambda_0\Sigma_1 v_1=  g_1 & \mbox{ in } \Omega, \\[6pt]
\dive(A_2 \nabla v_2) - i  \lambda_0\Sigma_2 v_2= g_2 & \mbox{ in } \Omega, 
\end{array}
\right.
\end{equation}
and  \eqref{lem1-v2-state1} and \eqref{lem1-v2-state2}  for some  $g = (g_1, g_2) \in [L^2(\Omega)]^2$ and $h \in H^{-1/2}(\Gamma)$, then 
\begin{equation*}
\| v\|_{\cH (\Omega)}^2 \le  C  {\cal  M}(v, g, h), 
\end{equation*}
for some positive constant  $C$ independent of $v$, $g $,  and $h$ where ${\cal M}(v, g, h)$ is defined in \eqref{def-M}. 
\end{lemma}

\noindent{\bf Proof.}  The proof is in the spirit of  the one of Lemma~\ref{lem-elementary-v2}. From \eqref{lem1-v2-state1}, we have, for all $0 < \gamma < 1$, 
\begin{equation}\label{coco-1}
 \lambda_0 \Big| \int_{\Omega} (\Sigma_2 - \Sigma_1) |v_{2}|^2 \Big| \le   \int_{\Omega} \gamma \lambda_0 |v_2 |^2 + C \gamma^{-1} \lambda_0 |w|^2 +  {\cal M}(v, g, h).
\end{equation}
Here and in what follows in the proof, $C$ denotes a positive constant depending only on $\Lambda_1$, $\Lambda_2$,  $\Omega$; it is independent of $\lambda_0$ and $\gamma$.  As in \eqref{pro-thm3-p3}, 
we deduce from \eqref{lem1-v2-state2} that 
\begin{equation}\label{coco-2}
\int_\Omega \langle A_1 \nabla w,  \nabla w \rangle + \int_{\Omega} \langle [A_{1} - A_{2}] \nabla v_{2} , \nabla v_{2} \rangle \le  \int_{\Omega} \gamma \lambda_0  |v_{2}|^2 + C \gamma^{-1} \lambda_0 |w|^2  + C {\cal M}(v, g, h).  
\end{equation}
Note that  $\int_{\Omega} |w|^2 \le C \int_{\Omega} |\nabla w|^2$ by  Poincar\'e's inequality, and 
\begin{equation}\label{coco-3}
\int_{\Omega} |v_2|^2 \le C \Big| \int_{\Omega} (\Sigma_2 - \Sigma_1) |v_{2}|^2 \Big| + C \int_{\Omega} d_\Gamma^\alpha |\nabla v_2|^2, 
\end{equation}
by Lemma~\ref{lem-decay3}. Fixing  $\gamma$ small enough in  \eqref{coco-1} and \eqref{coco-2}, adding \eqref{coco-1} and \eqref{coco-2}, and using \eqref{coco-3}, we derive that  if $\lambda_0$ is small enough, 
\begin{equation}\label{pro-thm3-p4}
\int_\Omega \langle A_{1} \nabla w,  \nabla w \rangle + \int_{\Omega} \langle [A_{1} - A_{2}] \nabla v_{2} , \nabla v_{2} \rangle + \lambda_0 \int_{\Omega} |v_2|^2
\le C {\cal M}(v, g, h).   
\end{equation}
The conclusion now follows from   \eqref{pro-thm3-p4} by noting that $w = v_1 - v_2$.  \proofend

\medskip 

We are ready to present

\medskip 
\noindent{\bf Proof of Proposition~\ref{pro-A1A2}}. The proof follows as in the one of Theorem~\ref{thm3} by considering the transformation 
\begin{eqnarray}
\begin{array}{rcccc}
T_2: & [L^2(\Omega)]^2 &\to&  [L^2(\Omega)]^2 \\[6pt]
& (f_1, f_2)  &  \mapsto & (u_1, u_2), 
\end{array}
\end{eqnarray}
where $(u_1, u_2) \in \cH(\Omega)$ is the unique viscosity solution of \eqref{Sys2}. The details are left to the reader. \proofend

\section{On the case $A_1 = A_2$ in a neighborhood of $\Gamma$ -  Proof of Theorem~\ref{thm2}} \label{sect-A1A2}

The condition $A_1 - A_2 \ge c d_\Gamma^\alpha I$ in a neighborhood of $\Gamma$ with $0 \le \alpha < 2$ plays a crucial role in establishing the compactness of $T$, more precisely $T_1$ and  $T_2$  (see Remarks~\ref{rem-com1} and \ref{rem-com2}). To be able to deal with case $A_1 = A_2$ in a neighborhood of $\Gamma$, we  make some modifications on $T$. The idea is to take into account the fact that  $u_1 - u_2 \in H^1_0(\Omega)$ which is more regular $u_1$ and $u_2$ which are in general not in $H^1(\Omega)$.
A modification on $T$ was also used in the work of Sylvester's \cite{Sylvester}.  Nevertheless, the modification in \cite{Sylvester} requires the condition $A_1 = A_2$ globally in $\Omega$  and  does not fit  for our situation.  

\medskip 
The motivation for reformulating the problem is as follows. Let $\lambda \in \mC$ be an eigenvalue of the ITE problem and let  $(u_1, u_2) \in [H^1(\Omega)]^2$ be a corresponding pair of eigenfunctions. Then 
 \begin{equation}\label{pro1a-1-1}  
 \left\{\begin{array}{cll}
 \dive(A_1 \nabla u_1) - \lambda\Sigma_1 u_1= 0 & \mbox{ in } \Omega, \\[6pt]
 \dive(A_2 \nabla u_2) - \lambda\Sigma_2 u_2= 0 & \mbox{ in } \Omega, \\[6pt]
  u_1 =u_2, \quad  A_1 \nabla u_1\cdot \nu = A_2 \nabla u_2\cdot \nu & \text{ on }\Gamma.
 \end{array} \right. 
 \end{equation}
Fix $\lambda_0 \neq 0$ and set $ w = u_1 - u_2$ in $\Omega$.  From \eqref{pro1a-1-1}, we have
\begin{align*}
\dive(A_1 \nabla w) - \lambda_0 \Sigma_1 w = &  (\lambda - \lambda_0) \Sigma_1 u_1  - \dive(A_1 \nabla u_2) + \lambda_0 \Sigma_1 u_2 \\[6pt]
 = &  (\lambda - \lambda_0) \Sigma_1 (u_1 - u_2)  -\lambda\Sigma_2u_2   -  \dive([A_1 - A_2]\nabla  u_2)+ \lambda \Sigma_1 u_2 \quad \mbox{ in } \Omega. 
\end{align*}
Define   \footnote{The goal is to eliminate $\Sigma_1 u_2$ from the equation of $\hat w$.}
$$
\hat w = w + \lambda u_2/ \lambda_0 \quad   \mbox{ in } \Omega. 
$$
Then 
\begin{equation}\label{motivation-1-v3}
\dive(A_1 \nabla \hat w) - \lambda_0 \Sigma_1 \hat w =   (\lambda - \lambda_0) \Sigma_1 (u_1 - u_2)    + \frac{\lambda }{\lambda_0}(\lambda-\lambda_0) \Sigma_2 u_2
 + \frac{ \lambda - \lambda_0}{ \lambda_0} \dive([A_1 - A_2]\nabla  u_2)    \mbox{ in } \Omega.  
\end{equation}
It is clear that 
\begin{equation}\label{motivation-2-v3}
\dive(A_2 \nabla u_2) - \lambda_0 \Sigma_2 u_2 = (\lambda - \lambda_0) \Sigma_2 u_2 \quad  \mbox{ in } \Omega. 
\end{equation}
Set 
\begin{equation*}
U_1 = \hat w = (u_1 - u_2) + \lambda u_2/ \lambda_0 \quad \mbox{ and } \quad U_2 =   \lambda u_2/ \lambda_0 \quad  \mbox{ in } \Omega. 
\end{equation*}
It follows from \eqref{motivation-1-v3} and \eqref{motivation-2-v3} that, in $\Omega$,  
\begin{multline}\label{pro1a-1}  
 \dive(A_1 \nabla U_1) - \lambda_0 \Sigma_1 U_1 \\[6pt]
 =   (\lambda - \lambda_0) \Sigma_1 (U_1 - U_2)   + (\lambda - \lambda_0)  \Sigma_2 U_2 + \frac{\lambda - \lambda_0}{ \lambda} \dive \big( [A_1 - A_2] \nabla U_2\big)
 \end{multline}
and
\begin{equation}\label{pro1a-2}  
 \dive(A_2 \nabla U_2) - \lambda_0 \Sigma_2 U_2= (\lambda - \lambda_0) \Sigma_2 U_2.
  \end{equation}
Since $A_1 = A_2$ near $\Gamma$, we also have 
 \begin{equation}\label{pro1b-1}
  U_1 =U_2, \quad  A_1 \nabla U_1\cdot \nu = A_2 \nabla U_2\cdot \nu ~\text{ on }\Gamma.
 \end{equation}
One can rewrite \eqref{pro1a-1}, \eqref{pro1a-2}, and \eqref{pro1b-1}  under the form  
$$
T_{3, \lambda} (U_1, U_2) = \frac{1}{\lambda -  \lambda_0} (U_1, U_2). 
$$
Here 
$$
T_{3, \lambda}(f_1, f_2) := (v_1, v_2), 
$$
where $(v_1, v_2)$ satisfies the system
\begin{equation}\label{pro-A}  
 \left\{\begin{array}{cll}
 \dive(A_1 \nabla v_1) - \lambda_0 \Sigma_1 v_1=   \dsp   \Sigma_1 (f_1 - f_2)  +   \Sigma_2 f_2 + \lambda^{-1} \dive \big( [A_1 - A_2] \nabla f_2 \big)& \mbox{ in } \Omega, \\[6pt]
 \dive(A_2 \nabla v_2) - \lambda_0 \Sigma_2 v_2 =  \Sigma_2 f_2  & \mbox{ in }\Omega, \\[6pt]
  v_1 = v_2, \quad  A_1 \nabla v_1\cdot \nu = A_2 \nabla v_2\cdot \nu & \text{ on }\Gamma.
\end{array}\right. 
 \end{equation}
The discreteness of $T$ can now be derived from  the discreteness of $T_{3, \lambda}$. 
To this end, we first  establish the well-posedness of the  system 
\begin{equation}\label{Sys3}  
 \left\{\begin{array}{cll}
 \dive(A_1 \nabla v_1) - \lambda_0 \Sigma_1 v_1=   g_1 + \dive (G_1)& \mbox{ in } \Omega, \\[6pt]
 \dive(A_2 \nabla v_2) - \lambda_0 \Sigma_2 v_2 =  g_2  & \mbox{ in }\Omega, \\[6pt]
  v_1 = v_2, \quad  A_1 \nabla v_1\cdot \nu = A_2 \nabla v_2\cdot \nu & \text{ on }\Gamma,
\end{array}\right. 
 \end{equation}
for appropriate functions $g_1, g_2$ and vector field $G_1$.  We follows the strategy used in the proof of Theorem~\ref{thm1}.  We first introduce some notations. Given $s \in \mR$, we denote 
$$
L^2(\Omega, d_\Gamma^s): = \Big\{\psi \in L^2_{\loc}(\Omega);  \| \psi \|_{L^2(\Omega, d_\Gamma^s)} < + \infty \Big\}, 
$$
where 
$$
\| \psi \|_{L^2(\Omega, d_\Gamma^s)}^2: = \int_{\Omega} d_\Gamma^s |\psi|^2.  
$$
We also define, for a given  $0 \le \beta_1 < 2$,  
$$
\cH_0(\Omega)  = \Big\{ (\psi_1, \psi_2) \in L^2_{\loc}(\Omega); \; \psi_1, \psi_2 \in L^2(\Omega, d_\Gamma^{\beta_1}) \mbox{ and } \psi_1 - \psi_2 \in L^2(\Omega, d_\Gamma^{-\beta_1})   \Big\}
$$
and 
$$
\| (\psi_1, \psi_2) \|_{\cH_0(\Omega)} : = \| (\psi_1, \psi_2) \|_{L^2(\Omega)} + \|\psi_1 - \psi_2 \|_{L^2(\Omega, d_\Gamma^{-\beta_1})}. 
$$

Here is a variant of Lemma~\ref{lem1}. 
\begin{lemma} \label{lem1-v3} Let $0\le  \beta_1 < 2$,  $\lambda_0 > 0$,  $g = (g_1, g_2) \in \cH_0(\Omega)$, and $G_1 \in [L^2(\Omega)]^d$ with  $\supp G_1 \subset \subset \Omega$.  Assume that   $v = (v_1, v_2)  \in [H^1(\Omega)]^2$  satisfy the system 
\begin{equation}\label{sys1-lem1-v3} 
\left\{ \begin{array}{cl}
\dive(A_1 \nabla v_1) - \lambda_0\Sigma_1 v_1=  g_1 + \dive (G_1) & \mbox{ in } \Omega, \\[6pt]
\dive(A_2 \nabla v_2) - \lambda_0\Sigma_2 v_2= g_2   & \mbox{ in } \Omega, \\[6pt]
v_1 = v_2, \quad  A_1 \nabla v_1\cdot \nu =  A_2 \nabla v_2\cdot \nu  &  \text{ on }\Gamma.
\end{array}
\right.
\end{equation}
We have, with $w = v_1 - v_2$, 
\begin{multline}\label{A3-1}
\int_{\Omega} \langle A_1 \nabla w, \nabla w \rangle + \lambda_0 \Sigma_1 |w|^2 \\[6pt]
= -  \int_{\Omega} (g_1 - g_2) \bar w  + \int_{\Omega} \langle G_1,  \nabla \bar w \rangle +   \int_{\Omega} \lambda_0 (\Sigma_2 - \Sigma_1) v_2 \bar w +  \int_{\Omega} \langle [A_2 - A_1] \nabla v_2, \nabla w \rangle 
\end{multline}
and 
\begin{multline}\label{lem1-part2-0-v3}
\int_{\Omega} \langle [A_1 - A_2] \nabla v_2, \nabla w \rangle  + \langle [A_1 - A_2] \nabla v_2,  \nabla v_2 \rangle   +  \lambda_0 (\Sigma_1 - \Sigma_2) |v_2|^2  \\[6pt]
 = \int_{\Omega}  g_2 \bar w  - (\bar g_1 - \bar g_2)  v_2  +  \lambda_0 (\Sigma_2 - \Sigma_1) v_2 \bar  w + \int_{\Omega} \langle  \nabla v_2, G_1 \rangle. 
\end{multline}
As a consequence of \eqref{A3-1} and \eqref{lem1-part2-0-v3}, we obtain 
\begin{equation}\label{lem1-part1-v3}
\int_{\Omega} \langle A_1 \nabla w, \nabla w \rangle +  \lambda_0 \Sigma_1 |w|^2  \le 4 {\cal N}(v, g, G_1)  +   \int_{\Omega} \big| \langle [A_2 - A_1] \nabla v_2, \nabla w \rangle \big| + \lambda_0 \big|(\Sigma_2 - \Sigma_1) v_2 \bar w\big|
\end{equation}
and 
\begin{multline}\label{lem1-part2-v3}
\int_{\Omega} \langle [A_1 - A_2] \nabla v_2,  \nabla v_2  \rangle  +  \lambda_0 (\Sigma_1 - \Sigma_2) |v_2|^2  \\
 \le 4 {\cal N}(v, g, G_1)  +  \int_{\Omega} \big| \langle [A_2 - A_1] \nabla v_2,  \nabla w \rangle  \big| +  \lambda_0 \big|(\Sigma_2 - \Sigma_1)   v_2 \bar  w \big|,  
\end{multline}
where
\begin{equation}\label{def-N}
{\cal N}(v, g, G_1) =  \int_{\Omega} |g| |w| + |g_1 - g_2| |v| + |G_1| |\nabla v_2|. 
\end{equation}
\end{lemma}

\noindent{\bf Proof.} The proof of Lemma~\ref{lem1-v3} follows closely to the one of Lemma~\ref{lem1} by noting that  one still has \eqref{part1-uniqueness} and \eqref{lem1-part2-0} in which $g_1$ is replaced by $g_1 + \dive (G_1)$ and $h = 0$, i.e, \eqref{A3-1} and \eqref{lem1-part2-0-v3} holds. The details are left to the reader.  \proofend

\medskip
We next establish a variant of Lemma~\ref{lem-stability}.

\begin{lemma} \label{lem-stability-v3}  Let $0 \le \beta <  \beta_1 < 2$ and assume that, for some $c, \, \tau > 0$, 
\begin{equation}\label{cond-F-thm1-*-v3}
			A_1  \ge  A_2  \quad \mbox{ and }  \quad  \Sigma_1 - \Sigma_2 \geq c d_\Gamma^\beta \mbox{ in } \Omega_\tau. 
\end{equation}
There exists $\Lambda_0 > 1$ such that if $\lambda_0 \ge  \Lambda_0$, and $v = (v_1, v_2)  \in [H^1_{\loc}(\Omega)]^2$  satisfies the system 
 \begin{equation}\label{sys1-cor1-v3} 
\left\{ \begin{array}{cl}
\dive(A_1 \nabla v_1) - \lambda_0\Sigma_1 v_1=  g_1 + \dive (G_1) & \mbox{ in } \Omega, \\[6pt]
\dive(A_2 \nabla v_2) - \lambda_0\Sigma_2 v_2= g_2 & \mbox{ in } \Omega,
\end{array}
\right.
\end{equation}
and  \eqref{lem1-part1-v3} and \eqref{lem1-part2-v3} for some  $g = (g_1, g_2) \in \cH_0(\Omega)$, and $G_1 \in [L^2(\Omega)]^d$ with $\supp G_1 \subset \subset \Omega$,  
then, for $w = v_1 - v_2$ in $\Omega$,  
\begin{multline}\label{est-lem-stability-v3}
\int_{\Omega_\tau} \langle [A_1 - A_2] \nabla v_2, \nabla v_2 \rangle + \int_{\Omega_\tau} (\Sigma_1 - \Sigma_2) |v_2|^2 + \int_{\Omega \setminus \Omega_{\tau}}\big( |\nabla v_2|^2 + |v_2|^2 \big)
 + \int_{\Omega}  |\nabla w|^2 \\[6pt]    
\le C \Big( {\cal N}(v, g, G_1) + \| g_2\|_{L^2(\Omega \setminus \Omega_{\tau/2})}^2 \Big), 
\end{multline}
for some positive constant  $C$ independent of $v$, $g $, and $G_1$ where ${\cal N}(v, g, G_1)$ is defined by \eqref{def-N}. In particular, if $v \in [H^1(\Omega)]^2$ is a solution of \eqref{sys1-lem1-v3} then \eqref{est-lem-stability-v3} holds. 

 \end{lemma}

\noindent{\bf Proof.} The proof is in the spirit of  the one of Lemma~\ref{lem-stability} and  even simpler. Adding \eqref{lem1-part1-v3} and \eqref{lem1-part2-v3}, using the fact, for $\gamma >0$, 
\begin{equation*}
\gamma \int_{\Omega_\tau} \langle [A_1 - A_2] \nabla w, \nabla w \rangle + \frac{1}{\gamma} \int_{\Omega_\tau} \langle [A_1 - A_2] \nabla v_2,  \nabla v_2  \rangle  \ge 2 \int_{\Omega_\tau} \Big| \langle [A_1 - A_2] \nabla w, \nabla v_2 \rangle \Big|
\end{equation*}
and
\begin{equation*}
\gamma \int_{\Omega_\tau} (\Sigma_1 - \Sigma_2) |w|^2 + \frac{1}{\gamma} \int_{\Omega_\tau} (\Sigma_1 - \Sigma_2) |v_2|^2  \ge 2 \int_{\Omega_\tau}  (\Sigma_1 - \Sigma_2) |w v_2|, 
\end{equation*} 
and taking $\gamma > 1$ and close to 1 in the previous two inequalities, we obtain 
\begin{multline}
\int_{\Omega}|\nabla w|^2 +  \lambda_0 |w|^2 + \lambda_0 (\Sigma_1 - \Sigma_2) |v_2|^2  + \int_{\Omega_\tau} \langle [A_1 - A_2] \nabla v_2, \nabla v_2 \rangle \\[6pt]
\le C {\cal N}(v, g, G_1)  +  C \int_{\Omega \setminus \Omega_\tau} |\nabla v_2|^2 + \lambda_0 |v_2 w|,
\end{multline}
for some positive constant  $C$  independent of $v$, $g$, $G_1$,  and $\lambda_0$. 
Since, for $\gamma > 0$,  
\begin{equation}\label{eq-eq}
2 |v_2 w| \le \gamma |v_2|^2 + \gamma^{-1} |w|^2, 
\end{equation}
by applying Lemma~\ref{lem-decay} to $v_2$ in $\Omega \setminus \Omega_{\tau/2}$ after fixing  $\gamma > 1$ large enough in \eqref{eq-eq}, we have, for large $\lambda_0$, 
\begin{multline*}
\int_{\Omega_\tau} \langle [A_1 - A_2] \nabla v_2, \nabla v_2 \rangle + \int_{\Omega_\tau} (\Sigma_1  - \Sigma_2) |v_2|^2 + \int_{\Omega \setminus \Omega_{\tau}} \big(|\nabla v_2|^2 + |v_2|^2 \big) 
 + \int_{\Omega} |\nabla w|^2 \\[6pt]
\le C \Big( {\cal N}(v, g, G_1) + \| g_2\|_{L^2(\Omega \setminus \Omega_{\tau/2})}^2 \Big),  
\end{multline*}
for some positive constant $C$ independent of $v$, $g$, and $G_1$; 
which is the conclusion. \proofend

\medskip 

We now introduce functional spaces. In what follows in this section, we {\bf assume} that \eqref{cond1-F-thm2} holds in $\Omega_\tau$ for some $\tau > 0$ and {\bf set} $\beta_1 = (2 + \beta)/2$. 
Define  
\begin{equation}
\hat \cH_1 (\Omega): = \Big\{ (v_1, v_2) \in [H^1_{\loc}(\Omega)]^2; \;   v_1 - v_2 \in H^1_0(\Omega), \;   \| (v_1, v_2) \|_{\hat \cH_1(\Omega)} < + \infty \Big\}, 
\end{equation}
where the norm is generated from the following scalar product, with $u = (u_1, u_2)$ and $v = (v_1, v_2)$,
\begin{equation}\label{scalar-hat-H1}
\langle u, v \rangle_{\hat \cH_1(\Omega)} =     \int_{\Omega} \nabla (u_1 - u_2) \nabla (\bar v_1 - \bar v_2) + \int_{\Omega} d_\Gamma^{\beta_1 + 2} \nabla u_2 \nabla \bar v_2 + \int_{\Omega \setminus \Omega_\tau } u \bar v + 
\int_{\Omega_\tau}  (\Sigma_1 - \Sigma_2) u \bar v  
\end{equation}
One can verify that $\hat \cH_1(\Omega)$ is a  Hilbert space.

\medskip 
We are ready to give the 
\begin{definition}\label{def-weak-3} Let $(g_1, g_2) \in  \cH_0(\Omega)$ and $G_1 \in [L^2(\Omega)]^d$ such that $\supp G_1 \subset \Omega \setminus \Omega_\tau$. A pair of functions $(v_1, v_2) \in \hat \cH_1(\Omega)$ is called a weak solution of \eqref{Sys3} if 
\begin{equation}\label{sys-def-weak-3}
\left\{ \begin{array}{cl}
\dive(A_1 \nabla v_1) - \lambda_0\Sigma_1 v_1=  g_1 + \dive (G_1) & \mbox{ in } \Omega, \\[6pt]
\dive(A_2 \nabla v_2) - \lambda_0\Sigma_2 v_2= g_2 & \mbox{ in } \Omega, \\[6pt]
(A_1 \nabla v_1-  A_2 \nabla v_2) \cdot \nu = 0  &  \text{ on }\Gamma. 
\end{array}
\right.
\end{equation}
\end{definition}

\begin{remark} \fontfamily{m} \selectfont
Since $\dive (A_1 \nabla v_1-  A_2 \nabla v_2) \in L^2 (\Omega \setminus \Omega_\tau)$ and $A_1 \nabla v_1-  A_2 \nabla v_2 = A_1 \nabla (v_1 - v_2) + (A_1 - A_2) \nabla v_2 \in L^2(\Omega)$ by \eqref{scalar-hat-H1}. The last identity on $\Gamma$ in \eqref{sys-def-weak-3} makes sense. 
\end{remark}

Applying Lemma~\ref{lem-stability-v3}, we can prove

\begin{lemma} \label{lem-elementary-v3}  Let $g= (g_1, g_2) \in  \cH_0(\Omega)$ and $G_1 \in L^2(\Omega)$ with $\supp G_1 \subset  \Omega \setminus \Omega_\tau$. 
There exists $\Lambda_0 > 1$ such that if $\lambda_0 \ge  \Lambda_0$, then for all $0< \delta < 1$, there exists a unique solution $v_\delta = (v_{1, \delta}, v_{2, \delta})  \in [H^1(\Omega)]^2$  satisfy the system 
 \begin{equation}\label{sys1-cor1-v3-1} 
\left\{ \begin{array}{cl}
\dive\big( (1 + \delta)A_1 \nabla v_{1, \delta} \big) - \lambda_0 (1 + \delta) \Sigma_1 v_{1, \delta}=  g_1 + \dive (G_1) & \mbox{ in } \Omega, \\[6pt]
\dive(A_2 \nabla v_{2, \delta}) - \lambda_0\Sigma_2 v_{2, \delta}= g_2 & \mbox{ in } \Omega,\\[6pt]
v_{1, \delta } = v_{2, \delta}, \quad (1+ \delta) A_1 \nabla v_{1, \delta} \cdot \nu = A_2 \nabla v_{2, \delta} \cdot  \nu & \mbox{ on } \Gamma. 
\end{array}
\right.
\end{equation}
Moreover,
\begin{equation}\label{est-lem-stability-v3-1}
\int_{\Omega} \delta \big( |\nabla v_{2, \delta}|^2 + |v_{2, \delta}|^2 \big) + \| v_\delta\|_{\hat \cH_1(\Omega)}^2
\le C \Big( \| g\|_{ \cH_0(\Omega)}^2 + \| G_1\|_{L^2(\Omega)}^2 \Big), 
\end{equation}
for some positive constant  $C$ independent of $\delta$, $g$, and $G_1$. As a consequence, there exists a weak solution $v = (v_1, v_2) \in \hat \cH_1(\Omega)$ of \eqref{Sys3} such that 
\begin{equation}\label{part2-pro-E1-v3}
\|v \|_{\hat \cH_1(\Omega)} \le C \Big( \| g \|_{ \cH_0 (\Omega)} + \|G_1 \|_{L^2(\Omega )}^2 \Big). 
\end{equation}

\end{lemma}

\noindent{\bf Proof.} Assume that $v_\delta = (v_{1, \delta}, v_{2, \delta}) \in [H^1(\Omega)]^2$ is a solution of \eqref{sys1-cor1-v3-1}. Set $w_\delta = v_{1, \delta} - v_{2, \delta}$.  Applying Lemma \ref{lem-stability-v3} to $(A_1, A_2) = \big((1 + \delta) A_1, A_2 \big)$ and  $(\Sigma_1, \Sigma_2) = \big( (1 + \delta) \Sigma_1, \Sigma_2 \big)$,  one obtains, for large $\lambda_0$, \begin{multline}\label{N-1}
\int_{\Omega_\tau} \langle \big[(1 + \delta)A_1 - A_2 \big] \nabla v_{2, \delta}, \nabla v_{2, \delta} \rangle + \int_{\Omega_\tau} \big((1 + \delta)\Sigma_1  - \Sigma_2\big) |v_{2, \delta}|^2  \\[6pt]  + \int_{\Omega \setminus \Omega_{\tau}} \big(|\nabla v_{2, \delta}|^2 + |v_{2, \delta}|^2 \big) 
+ \int_{\Omega} |\nabla w_\delta|^2 
\le C \Big( {\cal N}(v_\delta, g, G_1) + \| g_2\|_{L^2(\Omega \setminus \Omega_{\tau/2})}^2 \Big). 
\end{multline}
We have
$$
\int_{\Omega} |g| |w_{\delta}| \le  \Big(\int_{\Omega}d_{\Gamma}^{\beta_1}  |g|^2   \Big)^{1/2} \Big(\int_{\Omega} d_{\Gamma}^{-\beta_1} |w_\delta|^{2}   \Big)^{1/2}, 
$$
$$
\int_{\Omega} |g_1 - g_2| |v_\delta| \le  \Big(\int_{\Omega} d_{\Gamma}^{-\beta_1} |g_1 - g_2|^2   \Big)^{1/2} \Big(\int_{\Omega} d_{\Gamma}^{\beta_1}  |v_\delta|^2  \Big)^{1/2}, 
$$
$$
\int_{\Omega} |G_1| |\nabla v_{2, \delta}| \le \Big( \int_{\Omega \setminus \Omega_\tau} |G_1|^2 \Big)^{1/2} \Big( \int_{\Omega \setminus \Omega_\tau} |\nabla v_{2, \delta} |^2 \Big)^{1/2}. 
$$
Noting that  $w_\delta \in H^1_0(\Omega)$ and hence, by the Hardy inequality, we obtain 
\begin{equation*}
\int_{\Omega} d_\Gamma^{-2}|w_{\delta}|^2 \le C_\Omega \int_{\Omega} |\nabla w_\delta|^2.  
\end{equation*}
Since $0 \le \beta < \beta_1 < 2$, it follows from the definition of ${\cal N}$ that 
\begin{equation}\label{N-2}
{\cal N}(v_\delta, g, G_1) \le C \Big(  \|g \|_{\cH_0(\Omega)} + \| G_1\|_{L^2(\Omega \setminus \Omega_\tau)} \Big)\| v_\delta\|_{\hat \cH_1(\Omega)} 
\end{equation}
By Lemma~\ref{lem-multiplier}, we have  
\begin{equation}\label{N-4}
\int_{\Omega} d_\Gamma^{\beta + 2} |\nabla v_2|^2 \le C \int_{\Omega} d_\Gamma^\beta |v_2|^2 + d_\Gamma^{\beta} |g_2|^2.  
\end{equation}
Since $\beta_1 > \beta$, we derive from \eqref{N-4} that 
\begin{multline}\label{N-3}
\int_{\Omega_\tau} \langle \big[(1 + \delta)A_1 - A_2 \big] \nabla v_{2, \delta}, \nabla v_{2, \delta} \rangle + \int_{\Omega_\tau} \big((1 + \delta)\Sigma_1  - \Sigma_2\big) |v_{2, \delta}|^2 + \int_{\Omega \setminus \Omega_{\tau}} \big(|\nabla v_{2, \delta}|^2 + |v_{2, \delta}|^2 \big)   \\[6pt]
 + \int_{\Omega} |\nabla w|^2 + \| g\|_{\cH_0(\Omega)}^2 \ge C \Big( \| v_\delta\|_{\hat \cH_1(\Omega)}^2 +  \int_{\Omega} \delta \big( |\nabla v_{2, \delta}|^2 + |v_{2, \delta}|^2 \big) \Big)
\end{multline}
Combining \eqref{N-1}, \eqref{N-2}, and \eqref{N-3} yields \eqref{est-lem-stability-v3-1}. This in turn  implies the uniqueness of $v_\delta$. The existence of $v_\delta$ can be established via an approximation argument by first assuming that $g \in [L^2(\Omega)]^2$  and $G_1 \in [L^2(\Omega)]^d$ and then considering the general case; the existence in first case follows from Proposition~\ref{pro-existence1} \footnote{In Proposition~\ref{pro-existence1}, $G_1 = 0$; nevertheless, the same proof gives the same result in the case $G_1 \in L^2(\Omega)$ with $\supp G_1 \subset \subset \Omega$.}. The details of this fact are left to the reader. The existence and estimate  of a weak solution  $v \in \hat \cH_1(\Omega)$ of \eqref{Sys3} follows directly from the existence and the estimate of $v_\delta$. The details are omitted.\proofend

\medskip As in the proof of Theorem~\ref{thm1} in Section~\ref{sect-goodsign}, we  introduce the concept of viscosity solutions. 
 
\begin{definition}\label{def-WeakSolution3} Let $(g_1, g_2) \in  \cH_0(\Omega)$ and $G_1 \in L^2(\Omega)$ with $\supp G_1 \subset \Omega \setminus \Omega_\tau$. A weak solution  $(v_1, v_2) \in \hat \cH_1(\Omega)$ of \eqref{Sys3}  
is called a viscosity  solution if for any  $(\hv_1, \hv_2) \in \hat \cH_1(\Omega)$, $(\hg_1, \hg_2) \in  \cH_0(\Omega)$, and $\hG_1 \in [L^2(\Omega)]^d$ with $\supp \hG_1 \subset \Omega \setminus \Omega_\tau$ such that 
\begin{equation}\label{sys1-Uniqueness-0-v3} 
\left\{ \begin{array}{cl}
\dive(A_1 \nabla \hv_1) -  \lambda_0 \Sigma_1 \hv_1=  \hg_1 + \dive (\hG_1) & \mbox{ in } \Omega, \\[6pt]
\dive(A_2 \nabla \hv_2) -  \lambda_0\Sigma_2 \hv_2= \hg_2 & \mbox{ in } \Omega, \\[6pt]
A_1 \nabla \hv_1\cdot \nu - A_2 \nabla \hv_2\cdot \nu  = 0 &  \text{ on }\Gamma.
\end{array}
\right.
\end{equation}
we have, with $w = v_1 - v_2$ and $\hw = \hv_1 -\hv_2$,
\begin{multline}\label{lem1-part2-0-0-1-v3}
\int_{\Omega} \langle A_1 \nabla w, \nabla \hat w \rangle + \lambda_0 \Sigma_1 w \bar \hw = -  \int_{\Omega} (g_1 - g_2) \bar \hw  +   \int_{\Omega} \lambda_0 (\Sigma_2 - \Sigma_1) v_2 \bar \hw \\[6pt]
+  \int_{\Omega} \langle [A_2 - A_1] \nabla v_2, \nabla \hw \rangle + \int_{\Omega} G_1 \nabla \bar \hw
\end{multline}
and
\begin{multline}\label{lem1-part2-0-0-v3}
\int_{\Omega} \langle [A_1 - A_2] \nabla v_2, \nabla \hw \rangle  + \langle [A_1 - A_2] \nabla v_2,  \nabla \hv_2 \rangle   +  \lambda_0 (\Sigma_1 - \Sigma_2) v_2 \bar \hv_2  \\[6pt]
 = \int_{\Omega}  g_2 \bar \hw  - (\bar \hg_1 - \bar \hg_2)  v_2   + \bar \hG_1 \nabla v_2+  \lambda_0 (\Sigma_2 - \Sigma_1) v_2 \bar  \hw. 
\end{multline}
\end{definition}

We are now in a position to state and prove the  key result of this section. 
\begin{proposition}\label{pro-thm1-v3}  There exists $\Lambda_0 > 1$ such that for $\lambda_0 > \Lambda_0$ and  for  $g = (g_1, g_2) \in  \cH_0(\Omega)$ and $G_1 \in [L^2(\Omega)]^d$ with $\supp G_1 \subset \Omega \setminus \Omega_\tau$, there exists a unique viscosity solution $v =  (v_1, v_2) \in \hat \cH_1(\Omega)$ of \eqref{Sys3}. 
Moreover, 
\begin{equation}\label{est1-p-thm1-v3}
\| v \|_{\hat \cH_1 (\Omega)} \le C \Big( \| g \|_{ \cH_0(\Omega)} + \| G_1 \|_{L^2(\Omega)} \Big), 
\end{equation}
for some positive constant $C$ independent of $g$ and $G_1$. 
\end{proposition}

\noindent{\bf Proof.} We first prove that viscosity solutions exist. By Lemma~\ref{lem-elementary-v3}, there exists a unique solution  $v_{\delta} = (v_{1, \delta}, v_{2, \delta}) \in [H^1(\Omega)]^2$ ($0< \delta < 1$)  of the system
\begin{equation*}
\left\{ \begin{array}{cl}
\dive( (1 + \delta)A_1 \nabla v_{1, \delta}) - \lambda_0 (1 + \delta) \Sigma_1 v_{1, \delta}=  g_1 + \dive ( G_1 ) & \mbox{ in } \Omega, \\[6pt]
\dive(A_2 \nabla v_{2, \delta}) - \lambda_0\Sigma_2 v_{2, \delta} = g_2 & \mbox{ in } \Omega,\\[6pt]
v_{1, \delta } = v_{2, \delta}, \quad (1+ \delta) A_1 \nabla v_{1, \delta} \cdot \nu = A_2 \nabla  v_{2, \delta} \cdot \nu & \mbox{ on } \Gamma;  
\end{array}
\right.
\end{equation*}
moreover, 
\begin{equation}\label{est-lem-stability-v3-1-1}
\int_{\Omega} \delta \big( |\nabla v_{2, \delta}|^2 + |v_{2, \delta}|^2 \big) + \| v_\delta\|_{\hat \cH_1(\Omega)}^2
\le C \Big( \| g\|_{ \cH_0(\Omega)}^2 + \| G_1 \|_{L^2(\Omega)}^2 \Big), 
\end{equation}
for some positive constant  $C$ independent of $\delta$, $g $, and $G_1$.

 Set $w_\delta = v_{1, \delta} - v_{2, \delta}$. Let $v$ be the weak limit of $(v_{\delta_n})$ in $\hat \cH_1(\Omega)$ for some sequence $(\delta_n) \to 0$ such that $w_{\delta_n} \rightharpoonup w: = v_1 - v_2$ weakly in $H^1_0(\Omega)$. Then $v $ is a weak solution of \eqref{Sys3}.  We prove that $v$ is a viscosity solution. 
Multiplying the equation of $w_\delta$, 
\begin{multline*}
\dive\big((1 + \delta)A_1 \nabla w_\delta \big) - \lambda_0  (1 + \delta) \Sigma_1 w_\delta \\[6pt]
= g_1 + \dive (G_1) - g_2 -\lambda_0 \big(\Sigma_2 - (1 + \delta) \Sigma_1 \big) v_{2, \delta}  +  \dive ([A_2 -  (1 + \delta)A_1] \nabla v_{2, \delta}) \mbox{ in } \Omega, 
\end{multline*}
 by $\bar \hw$ and integrating in $\Omega$, we obtain 
\begin{multline}\label{coucou1}
\int_{\Omega} \langle (1 + \delta) A_1 \nabla w_\delta, \nabla \hat w \rangle + \lambda_0 (1 + \delta)\Sigma_1 w_\delta  \bar \hw 
= -  \int_{\Omega} (g_1 - g_2) \bar \hw  \\[6pt]+   \int_{\Omega} \lambda_0 \big(\Sigma_2 - (1 + \delta)\Sigma_1\big) v_{2, \delta} \bar \hw +  \int_{\Omega} \langle [A_2 - (1 + \delta)A_1 ] \nabla v_{2, \delta}, \nabla \hw \rangle + \int_{\Omega} G_1 \nabla \bar \hw. 
\end{multline}
 Similar to \eqref{lem1-part2-0-0}, we have 
\begin{multline}\label{coucou2}
\int_{\Omega} \langle [A_1 - A_2] \nabla v_{2, \delta}, \nabla \hw \rangle  + \langle [A_1 - A_2] \nabla v_{2, \delta},  \nabla \hv_2 \rangle   +  \lambda_0 (\Sigma_1 - \Sigma_2) v_{2, \delta} \bar \hv_2  \\[6pt]
 = \int_{\Omega}  g_2 \bar \hw  - (\bar \hg_1 - \bar \hg_2)  v_{2, \delta}   + \bar \hG_1 \nabla v_2+  \lambda_0 (\Sigma_2 - \Sigma_1) v_{2, \delta} \bar  \hw. 
\end{multline}
There is no term $(1 + \delta)$ in front of $A_1$ and $\Sigma_1$ in \eqref{coucou2} since we only use the equation of $v_2$ and the system of $\hat v$ here. Choosing  $\delta = \delta_n $, letting $n \to + \infty$ in \eqref{coucou1} and \eqref{coucou2} and using \eqref{est-lem-stability-v3-1-1}, we obtain  
 \eqref{lem1-part2-0-0-1-v3} and  \eqref{lem1-part2-0-0-v3}. Therefore,  $v$ is a viscosity solution.

 The uniqueness of viscosity solutions is now standard from its definition as in the proof of Proposition~\ref{pro-thm1}. \proofend
 
 \medskip 

Set
\begin{equation}
\hat \cH_0 (\Omega): = \Big\{ (v_1, v_2) \in [L^2_{\loc} (\Omega)]^2;  \| (v_1, v_2)\|_{\hat \cH_0(\Omega)}  < + \infty \Big\}, 
\end{equation}
where 
\begin{equation}\label{scalar-hat-H0}
\langle u, v \rangle_{\hat \cH_0(\Omega)} = \int_{\Omega \setminus \Omega_\tau} \nabla u_2 \nabla \bar v_2 + \int_{\Omega} d_\Gamma^{-\beta_1} (u_1 - u_2)(\bar  v_1 - \bar v_2) + \int_{\Omega} d_\Gamma^{\beta_1}u \bar v. 
\end{equation}
Fix $\lambda_0 >\Lambda_0$ where $\Lambda_0$ is the constant  in Proposition~\ref{pro-thm1-v3}. Define 
\begin{eqnarray}
\begin{array}{rcccc}
T_{3, \lambda}: & \hat \cH_0(\Omega) &\to& \hat \cH_0(\Omega) \\[6pt]
& (f_1, f_2)  &  \mapsto & (u_1, u_2), 
\end{array}
\end{eqnarray}
where $(u_1, u_2) \in \hat \cH_1(\Omega)$ is the unique viscosity solution of \eqref{Sys3} with 
\begin{equation}\label{notation-gG}
g_1  = \Sigma_1 (f_1 - f_2)  +   \Sigma_2 f_2, \quad G_1=  \lambda^{-1}  (A_1 - A_2) \nabla f_2,  \quad \mbox{ and } \quad g_2 = \Sigma_2 f_2.   
\end{equation}
We have 

\begin{lemma}\label{T3-compact} There exists $\Lambda_0 > 1$ such that if $\lambda_0 > \Lambda_0$ then $T_{3, \lambda}$ is a compact operator. 
\end{lemma}

\noindent{\bf Proof.}  Let $(f_{1, n}, f_{2, n})$ be an arbitrary bounded sequence in $\hat \cH_0(\Omega)$  and set 
$$
(u_{1, n}, u_{2, n})  = T_{3,\lambda}(f_{1, n}, f_{2, n}).
$$ 
Denote $g_{1, n}, g_{2, n}, G_{1, n}$ as in \eqref{notation-gG} where $(f_1, f_2)$ is replaced by $(f_{1, n}, f_{2, n})$.  
It follows from  \eqref{est1-p-thm1-v3} that  $(u_{1, n}, u_{2, n})$ is  bounded sequence in $\hat \cH_1(\Omega)$. Using the equations of $u_{1, n}$ and $u_{2, n}$, we derive that $(u_{1, n}, u_{2, n})$ is bounded in $H^1_{\loc}(\Omega)$. Since $\dive(A_2 \nabla u_{2, n}) - \lambda_0 \Sigma_2 u_{2, n} = g_{2, n}$, $(u_{2, n})$ is bounded in $H^1_{\loc}(\Omega)$,  and $(g_{2, n})$ is bounded in $L^2_{\loc}(\Omega)$, one might assume that 
\begin{equation}\label{compact-1}
(u_{2, n}) \mbox{ converges strongly in } H^1_{\loc}(\Omega).
\end{equation} 
By Hardy's inequality, we have, with $w_n = u_{1, n} - u_{2, n}$,  
\begin{equation*}
\int_{\Omega} d_\Gamma^{-2} |w_n|^2 \le C_\Omega \int_{\Omega} |\nabla w_n|^2. 
\end{equation*}
Since  the embedding $H^1 \subset L^2$ is compact,  applying Lemma~\ref{lem-compact-L2} below to $\psi_n = w_n$, $s = -2$, and  $t = - \beta_1$,  one might assume in addition that
\begin{equation}\label{compact-2}
(w_n) \mbox{ converges in } L^2(\Omega, d_{\Gamma}^{-\beta_1}). 
\end{equation}
Since  the embedding $H^1 \subset L^2$ is compact, applying  
Lemma~\ref{lem-compact-L2} below to $\psi_n = u_n$, $s = \beta$, and $t = \beta_1 > \beta$,  one might assume in addition that
\begin{equation}\label{compact-3}
(u_n) \mbox{ converges in } [L^2(\Omega, d_{\Gamma}^{\beta_1})]^2. 
\end{equation}
By Lemma~\ref{lem-multiplier}, we have 
\begin{equation*}
\int_{\Omega} d_\Gamma^{\beta + 2} |\nabla u_{2, n}|^2 \le C \int_{\Omega} d_\Gamma^\beta |u_{2, n}|^2 + d_\Gamma^{\beta} |g_{2, n}|^2 \le C.  
\end{equation*}
Using \eqref{compact-1} and applying  Lemma~\ref{lem-compact-L2} to $\psi_n = \nabla u_{2, n}$, $s = \beta + 2$, and $t = \beta_1 + 2$,  one can assume in addition that 
\begin{equation}\label{compact-4}
\nabla u_{2, n} \mbox{ converges in } [L^2(\Omega, d_{\Gamma}^{\beta_1 + 2})]^d. 
\end{equation}
The convergence of $(u_n)$ in $\hat \cH_0(\Omega)$ now follows from \eqref{compact-1}, \eqref{compact-2},  \eqref{compact-3}, and \eqref{compact-4}. \proofend

\begin{remark} \fontfamily{m} \selectfont  The embedding $\hat \cH_1(\Omega)$ into $\hat \cH_0(\Omega)$ is not compact. Nevertheless, $T_{3, \lambda}$ is compact as shown in Lemma~\ref{T3-compact}. 
\end{remark}

In the proof of Lemma~\ref{T3-compact}, we used  the following  compactness  result.

\begin{lemma}\label{lem-compact-L2} Let $m \ge 1$ and  $(\psi_n) \subset [L^2_{\loc}(\Omega)]^m$ be such that for every $K \subset \subset \Omega$, $(\psi_n |_{K})$ is relatively compact in $[L^2(K)]^m$. 
Assume that, for some $s \in \mR$, 
$$
\sup_{n } \| \psi_n \|_{[L^2(\Omega, d_\Gamma^s)]^m} < + \infty. 
$$
Then $(\psi_n)$ is relatively compact in $[L^2(\Omega, d_\Gamma^t)]^m$ for every $t>s$. 
\end{lemma}

\noindent{\bf Proof.}  The  proof is almost in the spirit of  the one \cite[Lemma 7]{Ng-WP}. For the sake of completeness, we present the proof. Set 
\begin{equation*}
C = \sup_{n } \| \psi_n \|_{[L^2(\Omega, d_\Gamma^s)]^m}^2. 
\end{equation*}
We have, for $\tau > 0$,  
\begin{equation}\label{L3-part2}
\int_{\Omega_\tau} d_\Gamma^t |\psi_n|^2 \le \tau^{t-s} \int_{\Omega_\tau} d_\Gamma^s |\psi_n|^2 \le C \tau^{t-2}. 
\end{equation}
Fix $\eps > 0$ arbitrary. Fix $\tau > 0$ small enough such that 
\begin{equation}\label{L3-part3}
\|\psi_n\|_{L^2(\Omega_\tau)} \le \eps/2 \quad \forall \, n \in \mN. 
\end{equation}
Such a constant  $\tau$ exists by \eqref{L3-part2}.  Since $(\psi_n)$ is relatively compact in $L^2(\Omega \setminus \Omega_\tau)$,  there exist $\psi_{n_1}, \cdots,  \psi_{n_k}$ such that 
\begin{equation}\label{L3-part4}
\Big\{ \psi_n \in L^2(\Omega \setminus \Omega_\tau);  \; n \in \mN \Big\} \subset \bigcup_{j=1}^k  \Big\{\psi \in L^2(\Omega \setminus \Omega_\tau) ; \| \psi- \psi_{n_j}\|_{L^2(\Omega \setminus \Omega_\tau, d_\Gamma^t)} \le \eps/2 \Big\}. 
\end{equation}
A combination of \eqref{L3-part3} and \eqref{L3-part4} yields 
\begin{equation*}
\big\{ \psi_n \in L^2(\Omega, d_\Gamma^t); \; n \in \mN \big\} \subset \bigcup_{j=1}^k \big\{\psi \in L^2(\Omega) ; \| \psi- \psi_{n_j}\|_{L^2(\Omega)} \le 3\eps/2 \big\}. 
\end{equation*}
Therefore,  $(\psi_n)$ is relatively  compact in $L^2(\Omega, d_\Gamma^t)$. \proofend

\medskip 
We are ready to give 

\medskip 

\noindent{\bf Proof of Theorem~\ref{thm2}.}  By  the spectral theory of a family of  compact analytic operators, the spectrum of $T_{3, \lambda}$ is discrete (see,  e.g., \cite[Theorem 8.92]{ReRo}). The conclusion follows by the definition of $T_{3, \lambda}$. \proofend

\section{On the case $A_1$ and $A_2$ satisfying the complementary condition on $\Gamma$ -  Proof of Theorem~\ref{thm4}} \label{sect-complementing}
	
This section is devoted to the proof of Theorem~\ref{thm4}. We first establish the well-posedness in $[H^1(\Omega)]^2$ of the following system (Proposition~\ref{pro-v4}), for $\lambda_0 > 0$ (large),  
 \begin{equation}\label{Sys4} 
 \left\{ \begin{array}{cl}
 \dive(A_1 \nabla u_1) - i  \lambda_0\Sigma_1 u_1= g_1 &\text{ in}~\Omega, \\[6pt]
 \dive(A_2 \nabla u_2) - i \lambda_0\Sigma_2 u_2= g_2 &\text{ in}~\Omega, \\[6pt]
  u_1 =u_2, \quad  A_1 \nabla u_1\cdot \nu = A_2 \nabla u_2\cdot \nu & \text{ on }\Gamma,
 \end{array}
 \right.
 \end{equation}
for a given pair $(g_1,g_2) \in [L^2(\Omega)]^2$ 
under the following three assumptions, which are  the assumptions of Theorem~\ref{thm4}:  

\begin{enumerate}

\item[$C1)$] $A_1, \; A_2, \;  \Sigma_1,  \;\Sigma_2$ are continous in a neighbourhood of $\Gamma$. 

\item[$C2)$] For all $x \in \Gamma$, $A_1(x)$ and $A_2(x)$ satisfy the following condition, with   $e = \nu(x)$,   
\begin{equation*}
\langle A_2(x) e, e \rangle \langle A_2(x) \xi, \xi \rangle  - \langle A_2(x) e,  \xi \rangle^2 \neq  \langle A_1(x) e, e \rangle \langle A_1(x) \xi, \xi \rangle  - \langle A_1(x) e, \xi \rangle^2 
\end{equation*}
for all $\xi \in {\cal P}(x) \setminus \{0 \}$, where
\begin{equation*}
{\cal P}(x) := \big\{\xi \in \mR^d; \langle \xi, e \rangle = 0  \big\}. 
\end{equation*}

\item[$C3)$]  For all $x \in \Gamma$, $A_1(x), \, A_2(x), \, \Sigma_1(x)$, and $ \Sigma_2(x)$ verify 
\begin{equation}\label{cond-complementary2}
\big\langle  A_1(x) \nu(x), \nu(x) \big\rangle \Sigma_1(x) \neq \big\langle  A_2(x) \nu(x), \nu(x) \big\rangle \Sigma_2(x). 
\end{equation}
 
\end{enumerate}
We then define the operator 
\begin{equation}\label{def-T}
T_4(f_1, f_2) = (u_1, u_2)  \mbox{ where $(u_1, u_2)$ is the unique solution of \eqref{Sys4} with $(g_1, g_2) = (\Sigma_1 f_1, \Sigma_2 f_2)$}.
\end{equation}
 
The basic ingredient in the proof of Proposition~\ref{pro-v4} is the following lemma:

\begin{lemma}\label{lem-Rd} Let $\lambda \ge 1$,  $A_1,A_2$ be two constant positive symmetric  matrices and let $\Sigma_1,\Sigma_2$ be two positive constants. Assume that 			 
\begin{equation}\label{cond-1-C-v4}
\langle A_2 e_d, e_d \rangle \langle A_2 \xi, \xi \rangle  - \langle A_2 e_d,  \xi \rangle^2 \neq  \langle A_1 e_d, e_d \rangle \langle A_1 \xi, \xi \rangle  - \langle A_1 e_d, \xi \rangle^2
\end{equation}
for all $ \xi \in P \setminus \{0 \}$,  where
\begin{equation*}
{\cal P} := \big\{\xi \in \mR^d; \langle \xi, e_d \rangle = 0  \big\}, 
\end{equation*}
and 
\begin{equation}\label{cond-2-C-v4}	
\big\langle  A_1 e_d, e_d \big\rangle \Sigma_1 \neq \big\langle  A_2 e_d, e_d \big\rangle \Sigma_2.
\end{equation}
Let $g = (g_1, g_2) \in [L^2(\mR^d_+)]^2$, $G =(G_1, G_2)  \in \big([L^2(\mR^d_+)]^d \big)^2$, $\varphi \in H^{1/2}(\mR^d_0)$, $\phi\in H^{-1/2}(\mR^d_0)$. Assume that  $v =(v_1,v_2) \in [H^1(\mR^d_+)]^2$ is a solution of the system 
\begin{equation}\label{pro3}
\left\{ \begin{array}{cl}
\dive(A_1 \nabla v_1) -i\lambda \Sigma_1 v_1=  g_1 +\dive(G_1)~~&\text{ in}~\mathbb{R}_+^d,\\[6pt]
\dive(A_2 \nabla v_2) -i\lambda \Sigma_2 v_2=  g_2 +\dive(G_2)~~&\text{ in}~\mathbb{R}_+^d,\\[6pt]
v_1-v_2= \varphi,~~ (A_1 \nabla v_1 - G_1)\cdot e_d - (A_2 \nabla v_2 - G_2)\cdot e_d =\phi & \text{ on } \mR^d_0.
\end{array} \right. 
\end{equation}
We have
\begin{multline}\label{est-key-v4}
\| v \|_{H^1(\mathbb{R}_+^d)}+ \lambda^{1/2}	\| v \|_{L^{2}(\mathbb{R}_+^d)}\\[6pt]
\leq C\left( \lambda^{-1/2} \|g\|_{L^{2}(\mathbb{R}_+^d)} + \|G\|_{L^{2}(\mathbb{R}_+^d)} + \| \varphi\|_{H^{1/2}(\mR^{d}_0)} + \lambda^{1/4} \| \varphi\|_{L^2(\mR^{d}_0)}
+  \|\phi \|_{H^{-1/2}(\mR^d_0)}\right), 
\end{multline}
where $C$ is a positive constant depending only on $d$,  $A_1$, $A_2$, $\Sigma_1$,  and $\Sigma_2$. 
\end{lemma}
	
Here and in what follows in this section, 
$$
\mR^d_+: = \big\{ x = (x', x_d) \in \mR^{d-1} \times \mR; \; x_d  > 0 \big\} \quad \mbox{ and }  \quad 
\mR^d_0: = \big\{ x = (x', x_d) \in \mR^{d-1} \times \mR; \; x_d  = 0 \big\}. 
$$

\noindent{\bf Proof.}   Let $u_j \in H^1(\mR^d_+)$ ($j=1, 2$) be the unique solution of  the system
\begin{equation*}
\dive(A_j \nabla u_j) -i\lambda \Sigma_j u_j=  g_j +\dive(G_j)  \quad  \mbox{ in } \mR^d_+, \quad \mbox{ and } \quad (A_j \nabla u_j - G_j ) \cdot e_d  = \phi_j, 
\end{equation*}
with $\phi_1 = \phi$ and $\phi_2 = 0$. 
Multiplying the equation of $u_j$ by $\bar u_j$ and integrating on $\mR^d_+$, we have 
\begin{equation*}
\int_{\mR^d_+} \langle A_j \nabla u_j, \nabla u_j \rangle + i \lambda \Sigma_j |u_j|^2 = \int_{\mR^d_+} G_j \nabla \bar u_j - g_j \bar u_j + \int_{\mR^d_0} \phi_j \bar u_j; 
\end{equation*}
which  implies 
\begin{equation}\label{toto1-C-v4}
\|\nabla u_j  \|_{L^{2}(\mathbb{R}_+^d)}+ \lambda^{1/2}	\| u_j \|_{L^{2}(\mathbb{R}_+^d)} \leq C\Big( \lambda^{-1/2} \|g_j\|_{L^{2}(\mathbb{R}_+^d)} + \|G_j \|_{L^{2}(\mathbb{R}_+^d)}  + \|  \phi_j\|_{H^{-1/2}(\mR^d_0)}\Big).  
\end{equation}
Since, for $j=1, 2$,  
$$
\int_{\mR^{d-1}} |u_j(x')|^2 \, dx' \le  2 \int_{\mR^d_+}  |\partial_{x_d} u_j | |u_j|,  
$$
it follows that, for $j=1, 2$, 
\begin{equation}\label{toto2-C-v4}
\| u_j \|_{L^2(\mR^d_0)}^2 \le 2 \| u_j \|_{L^2(\mR^d_+)} \| \nabla u_j\|_{L^2(\mR^d_+)}. 
\end{equation}
A combination of \eqref{toto1-C-v4} and \eqref{toto2-C-v4} yields, for $j=1, 2$,  
\begin{multline}\label{toto3-C-v4}
\|u_j  \|_{H^1(\mathbb{R}_+^d)}+ \lambda^{1/2}	\| u_j \|_{L^{2}(\mathbb{R}_+^d)}  + 
\| u_j \|_{H^{1/2}(\mR^d_0)} + \lambda^{1/4} \| u_j \|_{L^2(\mR^d_0)} \\[6pt]
\leq C\Big( \lambda^{-1/2} \|g_j\|_{L^{2}(\mathbb{R}_+^d)} + \|G_j \|_{L^{2}(\mathbb{R}_+^d)}  + \|  \phi_j\|_{H^{-1/2}}\Big).  
\end{multline}
By considering the system of $(v_1 - u_1, v_2 - u_2)$, one might assume that $g_1 = g_2 = 0$, $G_1 = G_2 = 0$,  and $\phi = 0$. In what follows, we make this assumption. 

\medskip 
Let $\hat v_j(\xi', t)$ for $j =1, 2$ and $\hat \varphi(\xi', t)$ be the Fourier transform of $v_j$ and $\varphi$ with respect to $x' \in \mR^{d-1}$, i.e., for $(\xi', t) \in \mR^{d-1} \times (0, + \infty)$, 
\begin{equation*}
\hat v_j(\xi', t) = \int_{\mR^{d-1}} v_j(x', t) e^{- i x' \cdot \xi' }\, dx' \quad  \mbox{ for } j=1, 2, \quad \mbox{ and } \quad  
\hat \varphi(\xi', t) = \int_{\mR^{d-1}} \varphi(x') e^{- i x' \cdot \xi' }\, dx'. 
\end{equation*}
Since 
\begin{equation*}
\dive(A_j \nabla v_j ) - i \lambda \Sigma_j v_j = 0 \mbox{ in } \mR^{d}_+,   
\end{equation*}
it follows that 
\begin{equation}\label{eq-hv}
a_j \hv_j''(t) + 2 i b_j \hv_j'(t) - (c_j + i \lambda \Sigma_j) \hv_j(t) = 0 \mbox{ for } t > 0,  
\end{equation}
where 
\begin{equation}\label{pro-aj-C-v4-0}
a_j  = (A_j)_{d,d}, \quad b_j =  \sum_{k=1}^{d-1} (A_j)_{d, k} \xi_k, \quad \mbox{ and } \quad c_j = \sum_{k=1}^{d-1} \sum_{l =1}^{d-1} (A_j)_{k, l} \xi_k \xi_l.  
\end{equation}
Here $(A_j)_{k, l}$ denotes the $(k, l)$ component of $A_j$ for $j=1, 2$ and the symmetry of $A_j$ is used. 
Define, for $j=1, 2$,  
\begin{equation}\label{def-Deltaj-v4}
\Delta_j = - b_j^2 +  a_j (c_j + i \lambda \Sigma_j). 
\end{equation}
Denote $\xi = (\xi', 0)$.  Since $A_j$ is symmetric and positive, it is clear that, for $j=1, 2$,  
\begin{equation}\label{est1-C-v4-0}
a_j = \langle A_j e_d, e_d \rangle > 0, \; \;  b_j = \langle A_j \xi, e_d \rangle,    \; \;  \mbox{ and }  \; \;  a_j c_j - b_j^2 =  \langle A_j e_d, e_d \rangle \langle A_j \xi, \xi \rangle - \langle A_j e_d, \xi \rangle^2 > 0.  
\end{equation}
For $j=1, 2$, let $\sqrt{\Delta_j}$ denote the square root of $\Delta_j$ with positive real part and set 
\begin{equation*}
\eta_j  =( - i b_j - \sqrt{\Delta_j})/ a_j. 
\end{equation*}
Since $\hat v_j(\xi', t) \in L^2(\mR^d_+)$, we derive from \eqref{eq-hv} that, for $j=1, 2$,  
\begin{equation}\label{def-vj-v4}
\hv_j(\xi', t) = \alpha_j (\xi') e^{\eta_j t}, 
\end{equation}
for some $\alpha_j  \in L^2(\mR^d_0)$.  Using the fact that $v_1 - v_2 = \varphi $ and $A_1 \nabla v_1 \cdot e_d - A_2 \nabla v_2 \cdot e_d = 0$ on $\mR^d_0$, we derive that 
\begin{equation}\label{equation-C-v4}
\alpha_1  (\xi' ) - \alpha_2 (\xi') = \hat \varphi(\xi')  \quad \mbox{ and } \quad \alpha_1 (\xi') \langle i A_1 \xi + \eta_1 A_1 e_d, e_d \rangle  - \alpha_2 (\xi') 
 \langle i A_2 \xi + \eta_2 A_2 e_d, e_d \rangle  = 0. 
\end{equation}
Since, by \eqref{est1-C-v4-0}, 
$$
\langle A_j \xi, e_d \rangle - \langle A_j e_d, e_d \rangle  b_j/ a_j = 0 \mbox{ for } j=1, 2, 
$$
the last identity of \eqref{equation-C-v4} implies 
\begin{equation*}
\alpha_1(\xi') \sqrt{\Delta_1} = \alpha_2(\xi') \sqrt{\Delta_2}. 
\end{equation*}
Combining this identity and the first one of \eqref{equation-C-v4} yields 
\begin{equation}\label{def-alphaj-v4}
\alpha_1(\xi') = \frac{\hat \varphi(\xi') \sqrt{\Delta_2}}{\sqrt{\Delta_2} - \sqrt{\Delta_1}}. 
\end{equation}
Note that, by \eqref{cond-1-C-v4},  \eqref{cond-2-C-v4}, \eqref{def-Deltaj-v4}, and  \eqref{est1-C-v4-0},
\begin{align*}
|\Delta_2 - \Delta_1|^2 &=  \Big(\langle A_2 e_d, e_d \rangle \langle A_2 \xi, \xi \rangle  - \langle A_2 e_d,  \xi \rangle^2 -  \langle A_1 e_d, e_d \rangle \langle A_1 \xi, \xi \rangle  + \langle A_1 e_d, \xi \rangle^2  \Big)^2  \\[6pt]
& \quad + \lambda^2 \Big( \big\langle  A_1 e_d, e_d \big\rangle \Sigma_1 - \big\langle  A_2 e_d, e_d \big\rangle \Sigma_2 \Big)^2 \\[6pt]
& \ge C (|\xi|^4 + \lambda^2)
\end{align*}
and
\begin{equation*}
|\Delta_j| \le C (|\xi|^2 + \lambda). 
\end{equation*}
This implies 
\begin{equation}\label{est1-C-v4}
\frac{\sqrt{\Delta_2}}{\big|\sqrt{\Delta_2} - \sqrt{\Delta_1} \big|} = \frac{\sqrt{\Delta_2}\big(\sqrt{\Delta_2} + \sqrt{\Delta_1}\big)}{\big|\Delta_2 - \Delta_1 \big|}  \le  C. 
\end{equation}
Since, by Parseval's theorem,    
$$
\int_{\mR^d_+} |v_1|^2 = C_d \int_0^\infty \int_{\mR^{d-1}} |\hat v_1(\xi', t)|^2 \, d\xi' \, dt , 
$$
it follows from \eqref{def-vj-v4}, \eqref{def-alphaj-v4}, and  \eqref{est1-C-v4} that 
\begin{equation}\label{est2-C-v4}
\int_{\mR^d_+} |v_1|^2 \le C \int_{\mR^{d-1}} \int_0^\infty    |\hat \varphi(\xi')|^2|e^{2 \eta_1 t}| \, dt \, d \xi'. 
\end{equation}
Using the fact 
\begin{equation*}
\int_0^\infty |e^{2 \eta_1 t}| \, dt  \le \frac{1}{2 |\Re (\eta_1)|} \quad \mbox{ and } \quad  
\frac{1}{|\Re (\eta_1)|} \le \frac{C}{|\xi'| + \sqrt{\lambda}}, 
\end{equation*}
we deduce from \eqref{est2-C-v4} that 
\begin{equation*}
\int_{\mR^d_+} |v_1|^2 \le  \int_{\mR^{d-1}}  \frac{C|\hat \varphi(\xi')|^2}{|\xi'| + \sqrt{\lambda} } \, d \xi' \le C \lambda^{-1/2} \| \varphi\|_{L^2(\mR^{d}_0)}^2. 
\end{equation*}
which yields the corresponding  estimate for $\| v_1 \|_{L^2(\mR^d_+)}$. The estimates follows similarly. The details are left to the reader. \proofend

\begin{remark}  \fontfamily{m} \selectfont  The computation of $\hat v_j$ in Lemma~\ref{lem-Rd} has roots from the proof of  \cite[Proposition 1]{Ng-WP}. 
\end{remark}

\begin{proposition}\label{pro-v4}  Assume C1), C2), and C3).
There exists a constant $\Lambda_0>0$ depending only on $A_1, \, A_2,\,  \Sigma_1, \, \Sigma_2$, and $\Omega$, such that  for $\lambda >\Lambda_0$ and  $ g = (g_1,g_2) \in [L^{2}(\Omega)]^2$,  there exists a unique solution $v = (v_1,v_2) \in [H^1(\Omega)]^2$ satisfy 
\begin{equation}\label{pro5}
\left\{ \begin{array}{cl}
\dive(A_1 \nabla v_1) -i\lambda \Sigma_1 v_1=  g_1 ~~&\text{ in}~\Omega,\\[6pt]
\dive(A_2 \nabla v_2) -i\lambda \Sigma_2 v_2=  g_2 ~~&\text{ in}~\Omega,\\[6pt]
v_1-v_2=0,~~	A_1 \nabla v_1\cdot \nu-A_2 \nabla v_2\cdot \nu=0 & \text{ on }\Gamma. 
\end{array}\right. 
\end{equation}
Moreover, 
\begin{equation}\label{est-C1}
\|v \|_{H^1(\Omega)}+\lambda^{1/2} \| v \|_{L^{2}(\Omega)}\leq C \lambda^{-1/2}||g||_{L^{2}(\Omega)}, 	
\end{equation}
for some positive constant  $C$  independent of $g$ and $\lambda$.
\end{proposition}
	
\noindent{\bf Proof.}  We first assume the existence of $v \in [H^1(\Omega)]^2$ and derive the estimate for $v$.  Assume that  $(A_1, A_2, \Sigma_1, \Sigma_2) \in C(\bar \Omega_{ \tau_0})$ for some $\tau_0 > 0$.  Using local charts, applying Lemma~\ref{lem-Rd} with $\varphi = 0$ and $\phi = 0$, and involving the standard freezing coefficient technique, we have
\begin{equation}\label{part1-Thm4}
\| v\|_{H^1(\Omega_\tau)} + \lambda^{1/2} \|v\|_{L^2(\Omega_\tau)} \le C \Big(  \| v \|_{H^1(\Omega_{2 \tau} \setminus \Omega_\tau)}  + \lambda^{-1/2} \| (g_1, g_2)\|_{L^2(\Omega)} \Big),  
\end{equation}
if $\tau>0$ is small enough and $\lambda > 1$.  Here and in what follows $C$ denotes a positive constant independent of $v$, $g$, and $\lambda$; $C$ depends on $(A_1, A_2, \Sigma_1, \Sigma_2)$, $\Omega$,  and $\tau$. Fix such a positive constant  $\tau$. 
Applying Lemma~\ref{lem-decay-2} to $v$, we have 
\begin{equation}\label{part2-Thm4}
\| v \|_{H^1(\Omega  \setminus \Omega_\tau)} \le c_1 \exp(-c_2 \sqrt{\lambda}) \| v\|_{L^2(\Omega_\tau)} + c_1 \| g \|_{L^2(\Omega)}, 
\end{equation}
for some positive $c_1$ and $c_2$ independent of $\lambda$, $v$, and $g$. 
Combining \eqref{part1-Thm4} and \eqref{part2-Thm4} yields, for large $\lambda$, 
\begin{equation*}
\| v\|_{H^1(\Omega)} + \lambda^{1/2} \| v \|_{L^2(\Omega)} \le C \ \| g \|_{L^2(\Omega)}; 
\end{equation*}
which is \eqref{est-C1}. The existence of $v \in [H^1(\Omega)]^2$ follows from the uniqueness via the Fredholm theory  by noting  
the well-posedness in $[H^1(\Omega)]^2$ of the system
\begin{equation*}
\left\{ \begin{array}{cl}
\dive(A_1 \nabla v_1) -i\lambda \Sigma_1 v_1=  g_1 ~~&\text{ in}~\Omega,\\[6pt]
\dive(A_2 \nabla v_2) + i\lambda \Sigma_2 v_2=  g_2 ~~&\text{ in}~\Omega,\\[6pt]
v_1-v_2=0,~~	A_1 \nabla v_1\cdot \nu-A_2 \nabla v_2\cdot \nu=0 & \text{ on }\Gamma, 
\end{array}\right. 
\end{equation*}
for $g \in [L^2(\Omega)]^2$ by Lax-Milgram's theory.  
The details are left to the reader. \proofend

	
\medskip 	

We are ready to present 
\medskip 

\noindent{\bf Proof of Theorem~\ref{thm4}.} Fix $\lambda > \Lambda_0$ where $\Lambda_0$ is the positive constant  in Proposition~\ref{pro-v4}. Define 
			\begin{eqnarray}
			\begin{array}{rcccc}
			T_4: & [L^2(\Omega)]^2 &\to&  [L^2(\Omega)]^2 \\[6pt]
			& (f_1, f_2)  &  \mapsto & (u_1, u_2), 
			\end{array}
			\end{eqnarray}
			where $(u_1, u_2) \in [H^1(\Omega)]^2$ is the unique solution of \eqref{Sys4} with $(g_1, g_2) = (\Sigma_1 f_1, \Sigma_2 f_2)$. Since $H^1(\Omega)\subset L^2(\Omega)$ is compact, so  $T_4$ is  compact. The conclusion of Theorem~\ref{thm4} follows.  \proofend

								 \providecommand{\bysame}{\leavevmode\hbox to3em{\hrulefill}\thinspace}
								 \providecommand{\MR}{\relax\ifhmode\unskip\space\fi MR }
								 \providecommand{\MRhref}[2]{%
								 	\href{http://www.ams.org/mathscinet-getitem?mr=#1}{#2}
								 }
								 \providecommand{\href}[2]{#2}

\end{document}